\documentclass[11pt,a4paper]{article}
\usepackage[T1]{fontenc}
\usepackage[utf8]{inputenc}
\usepackage{lmodern}
\usepackage[margin=23mm]{geometry}
\usepackage{amsmath,amssymb,amsthm,mathtools,mathrsfs}
\usepackage{enumitem}
\usepackage{booktabs}
\usepackage{array}
\usepackage{graphicx}
\usepackage{xcolor}
\usepackage{microtype}
\usepackage[font=small]{caption}
\setlist{topsep=4pt plus 1pt,itemsep=1.5pt,parsep=1pt}
\usepackage[hidelinks]{hyperref}
\graphicspath{{./}{figures/}}

\theoremstyle{plain}
\newtheorem{theorem}{Theorem}[section]
\newtheorem{proposition}[theorem]{Proposition}

\newtheorem{corollary}[theorem]{Corollary}

\theoremstyle{definition}
\newtheorem{definition}[theorem]{Definition}
\newtheorem{problem}[theorem]{Problem}

\theoremstyle{remark}
\newtheorem{remark}[theorem]{Remark}

\newcommand{\R}{\mathbb{R}}

\newcommand{\E}{\mathbb{E}}
\newcommand{\Prob}{\mathbb{P}}
\newcommand{\Id}{\mathrm{I}}
\DeclareMathOperator{\Cov}{Cov}

\DeclareMathOperator{\dist}{dist}
\DeclareMathOperator{\supp}{supp}

\DeclareMathOperator{\Law}{Law}
\DeclareMathOperator{\conv}{conv}
\DeclareMathOperator{\Bis}{Bis}
\DeclareMathOperator{\arccosh}{arccosh}
\DeclareMathOperator{\arcsinh}{arcsinh}
\newcommand{\dd}{\,\mathrm{d}}
\newcommand{\lap}{\Delta}
\newcommand{\Hess}{\nabla^{2}}

\newcommand{\sg}{\mathsf{s}} 

\newcommand{\one}{\mathbf 1}

\title{\textbf{Speciation in score-based diffusion as a partition problem\\[2pt]  }}
\author{Farid Bozorgnia\\ \small Department of Mathematics, New Uzbekistan University, Tashkent\\ \small\texttt{f.bozorgnia@newuu.uz}}
\date{}

\begin{document}
\maketitle
\begin{abstract}
We study \emph{speciation}---the commitment of reverse-diffusion trajectories to data classes---as a partition problem in $\R^d$ for a fixed data law, with explicit dimension-free constants. Under an exact score, time reversal identifies reverse-SDE committors with Bayesian posteriors: the resulting Bayes partition is a Voronoi or power diagram for atomic data, whereas probability-flow ODE basins carry the exact class masses. We prove four main results. For $K$ separated compactly supported clusters, an explicit noise window yields exactly $K$ nondegenerate modes, alongside a sharp diameter threshold for unimodality and two-sided clone-agreement bounds. For two convexly separated classes, the Bayes interface is a real-analytic graph at every noise level, converging to the support bisector as noise vanishes, with transition-layer probability given by a coarea formula. For three equally weighted non-collinear atoms, the ODE basins are optimal-transport cells if and only if the triangle is equilateral. For an approximate score stopped at positive noise (with the exact posterior below it), the Bayes partition changes only inside the transition layer, up to the integrated score error. Numerical experiments compute finite-difference committors, layer widths, and ODE basins for interlocking arcs and three atoms.

\end{abstract}
\medskip
\noindent\textbf{Keywords:} score-based diffusion models; reverse SDE; probability-flow ODE; committor functions; speciation; semi-discrete optimal transport; Gaussian scale space.

\noindent\textbf{MSC 2020:} 60J60, 49Q22, 60H10, 35Q84, 62H30.
\section{Introduction}\label{sec:intro}

Score-based diffusion models \cite{SohlDickstein2015,Ho2020,Song2021SDE} generate data by reversing a continuous noising process. Forward in time, an Ornstein--Uhlenbeck (OU) process gradually corrupts a data sample $X_0\sim p_0$ until its distribution approaches a standard Gaussian. Backward in time, new samples are generated from Gaussian noise by integrating either a stochastic differential equation (the \emph{reverse SDE}) or a deterministic flow (the \emph{probability-flow ODE}), both driven by the score $\nabla\log p_t$ of the noised density $p_t$.

Information about the generated sample does not emerge uniformly along a trajectory. At high noise, a state carries almost no information about its endpoint. As the noise decreases, the trajectory passes through a scale where it commits to a class, cluster, or coarse semantic feature. This crossover is called \emph{speciation} and is often compared to symmetry breaking in statistical physics \cite{Raya2023,Biroli2024Regimes}. In fixed dimension, however, speciation is smooth rather than sharp, and different ways of measuring it need not occur at the same noise level (Remark~\ref{rem:three}). At a much smaller noise scale, a sampler driven by the empirical score of $n$ training points commits to a single training point. 

Knowing \emph{when} and \emph{where} trajectories commit matters in both practice and theory. In conditional generation, guidance \cite{Dhariwal2021} is effective mainly on an intermediate noise interval---it can be harmful at high noise and is largely unnecessary at low noise \cite{Kynkaanniemi2024}---making the speciation window a natural candidate for this interval. Beyond generative modelling, forward noising defines a Gaussian scale space whose modes merge as noise increases (monotonically in one dimension), linking speciation to scale-space theory, mean-shift clustering \cite{Lindeberg1994,FukunagaHostetler1975,Cheng1995}, and Silverman's multimodality test \cite{Silverman1981}.

\emph{Related work.} In statistical physics, Raya and Ambrogioni \cite{Raya2023} linked speciation to spontaneous symmetry breaking, while Biroli and M\'ezard \cite{BiroliMezard2023} and Biroli et al.\ \cite{Biroli2024Regimes} studied high-dimensional speciation via the probability that two independent continuations (\emph{clones}) from a shared noisy state reach the same class, connecting the speciation time to the data covariance. On the mathematical side, critical windows for mixtures have been bounded using forward total-variation estimates \cite{LiChen2024} and stochastic localization \cite{LiKaranChen2025}. For the deterministic sampler, the probability-flow map is generally not an optimal-transport map \cite{LavenantSantambrogio2022}, and its partition cells for discrete targets are studied in \cite{Pierret2026}, which constructs a planar four-point non-Laguerre example.

Our paper complements these works. Rather than taking a high-dimensional limit ($d\to\infty$), we fix an arbitrary compactly supported data law $p_0$ in $\R^d$ (possibly singular) and prove explicit bounds whose constants do not depend on the ambient dimension $d$---only on cluster separations, diameters, class weights, and (in Theorem~\ref{thm:interface}(d)) a lower mass bound. Instead of locating only a critical time, we describe both the modes of $p_t$ and the full geometry of the decision interface, and for the deterministic sampler we give a complete three-atom characterization at every positive noise level.

\emph{The partition viewpoint.} Our starting point is that speciation is a partition problem: at every forward time $t>0$, the reverse dynamics partition $\R^d$ by trajectory endpoints---softly through hitting probabilities (\emph{committors}) and sharply through the regions where one class is most probable. Under the exact score, a classical consequence of time reversal (Theorem~\ref{thm:committor}) identifies the reverse-SDE committor with the Bayesian posterior $\Prob(X_0\in G\mid X_t=x)$. As a result, the committor---a priori the solution of a singular backward Kolmogorov equation---is explicit, the stochastic decision boundaries are Bayes interfaces, and a class with posterior probability $1-\varepsilon$ remains the most probable one along the entire trajectory with probability at least $1-2\varepsilon$ (Corollary~\ref{cor:explicit}). For atomic data, this Bayes partition is simply the Voronoi diagram of the rescaled atoms (or a power diagram for unequal weights) at \emph{every} $t>0$.  Genuine free boundaries arise only in the zero-noise limit, for approximate scores, and in the frozen-time problem (Remark~\ref{rem:FB}).

\emph{Main results.} Within this framework, we prove four main results:
\begin{enumerate}[label=(\roman*),nosep,leftmargin=2.2em]
\item \emph{An explicit $K$-mode window.} For $K$ separated clusters of radius $r$ and pairwise separation $\Delta$, the noised density $p_t$ has exactly $K$ nondegenerate modes (one per cluster) for all noise levels $h_t\in[2r,h_+]$, and a single mode once $h_t$ exceeds half the support diameter; clone agreement is bounded from both sides at the same scales (Theorem~\ref{thm:Kmodes}, Corollary~\ref{cor:window}, Propositions~\ref{prop:logconcave} and~\ref{prop:clonewindow}).
\item \emph{Global geometry of the interface.} For two convexly separated classes, the Bayes interface is a real-analytic graph over a hyperplane at every $t>0$. Under a lower mass bound, it converges locally to the support bisector at rate $h_t^2\log(1/h_t)$ (or $h_t^2$ for atomic classes), and the transition-layer probability is given by a coarea formula (Theorems~\ref{thm:interface} and~\ref{thm:width}).
\item \emph{ODE basins and optimal transport.} The probability-flow ODE basins carry the exact class masses (Proposition~\ref{prop:odemass}). For two atoms they coincide with optimal-transport cells, with an interface parallel to the Bayes interface at an explicit distance (Proposition~\ref{prop:odesde}); for three equally weighted non-collinear atoms they are optimal-transport cells if and only if the triangle is equilateral (Theorem~\ref{thm:threeatoms}).
 
\end{enumerate}

 The organization of the paper is as follows.  Section~\ref{sec:setting} fixes the setting. Sections~\ref{sec:exact} and~\ref{subsec:twoatoms}--\ref{subsec:modes} assemble the posterior framework and benchmark examples (where time reversal, Tweedie's formula, the Robertson--Fryer criterion, one-dimensional mode monotonicity, and Duistermaat's non-monotone example in Proposition~\ref{prop:nonmonotone} are classical). The main new results are Theorem~\ref{thm:threeatoms} and Sections~\ref{subsec:Kmodes} and~\ref{sec:interfaces}; Section~\ref{sec:numerics} reports numerical experiments.

\section{Setting and problem formulation}\label{sec:setting}
\paragraph{Notation.} $|\cdot|$ is the Euclidean norm, $B(x,r)$ the open ball with centre $x$ and radius $r$, $B_R:=B(0,R)$, and $\bar B_R$ its closure. $\Id=\Id_d$ is the identity matrix, and $A\preceq B$ means that $B-A$ is positive semidefinite. $N(m,\Sigma)$ is the Gaussian law with mean $m$ and covariance $\Sigma$, $\Law(X)$ the law of a random variable $X$, $\mu*\nu$ the convolution of two measures, and $\mathcal H^{k}$ the $k$-dimensional Hausdorff measure.  $\Phi$ denotes the standard normal distribution function, $\varphi$ its density, and $s(z):=(1+e^{-z})^{-1}$ the logistic function. For $0<s\le t$, $\Phi_{t\to s}(x)$ denotes the state at forward time $s$ along the probability-flow ODE trajectory started at $x$ at forward time $t$ (Section~\ref{subsec:forward}); the flow map always carries its time subscript.  
\subsection{Forward process, posterior identities, and reverse dynamics}\label{subsec:forward}

\emph{Forward (noising) process.} Throughout, $p_0$ denotes the data distribution: a probability measure on $\R^d$ with compact support, which need not have a density (Section~\ref{subsec:data}). A data sample $X_0\sim p_0$ is noised by the Ornstein--Uhlenbeck process, in the normalization of \cite{Biroli2024Regimes}:
\begin{equation}\label{eq:OU}
\dd X_t=-X_t\dd t+\sqrt2\dd B_t,\qquad X_t\overset{d}{=}e^{-t}X_0+\sigma_tZ,\qquad \sigma_t^2:=1-e^{-2t},
\end{equation}
where $B$ is a Brownian motion in $\R^d$ and $Z\sim N(0,\Id)$, both independent of $X_0$. Thus $X_t$ is the data sample after noising for a time $t$: the signal is damped by $e^{-t}$ and Gaussian noise of variance $\sigma_t^2$ is added, so $X_t$ interpolates between the data ($t=0$) and pure noise $N(0,\Id)$ ($t\to\infty$). Indeed, the linear equation has the solution $X_t=e^{-t}X_0+\sqrt2\int_0^te^{-(t-s)}\dd B_s$, and the stochastic integral is centred Gaussian, independent of $X_0$, with covariance $2\int_0^te^{-2(t-s)}\dd s\,\Id=\sigma_t^2\Id$. Consequently, for every $t>0$ the law of $X_t$ has the smooth and strictly positive density
\[
p_t(x)=\int_{\R^d}(2\pi\sigma_t^2)^{-d/2}\exp\Big(-\frac{|x-e^{-t}z|^2}{2\sigma_t^2}\Big)\,p_0(\dd z),
\]
even when $p_0$ itself is singular. We call $p_t$ the \emph{noised density} and $s_t:=\nabla\log p_t$ its \emph{score}.

\emph{Noise level and data coordinates.} Multiplying \eqref{eq:OU} by $e^t$ removes the damping: $e^tX_t\overset{d}{=}X_0+h_tZ$, where
\begin{equation}\label{eq:h}
h_t:=e^t\sigma_t=\sqrt{e^{2t}-1}
\end{equation}
is the \emph{noise level}, the standard deviation of the noise in the units of the data. It increases from $0$ to $\infty$ with $t$, and $t=\tfrac12\log(1+h_t^2)=\tfrac12h_t^2+O(h_t^4)$ as $t\downarrow0$. We call $y=e^tx$ \emph{data coordinates}: the rescaled state $e^tX_t$ is the data sample observed through a Gaussian channel of noise level $h_t$, and its density is $p^{h_t}$, where $p^h:=p_0*N(0,h^2\Id)$. Since $X_t=e^{-t}(e^tX_t)$, the change of variables $x=e^{-t}y$ gives
\begin{equation}\label{eq:scaling}
p_t(x)=e^{td}\,p^{h_t}(e^tx).
\end{equation}
We pass freely between OU coordinates $(x,t)$ and data coordinates $(y,h)$ through $y=e^tx$ and $h=h_t$.

\emph{Posterior and Tweedie's identities.} By \eqref{eq:OU}, conditionally on $X_0=z$ the state $X_t$ has law $N(e^{-t}z,\sigma_t^2\Id)$. Bayes' rule therefore gives the law of the data sample $X_0$ given the noisy observation $X_t=x$:
\begin{equation}\label{eq:posterior}
\Pi_t(\dd z\mid x)=\frac{e^{-|x-e^{-t}z|^2/(2\sigma_t^2)}\,p_0(\dd z)}{\int e^{-|x-e^{-t}z'|^2/(2\sigma_t^2)}\,p_0(\dd z')}
=\frac{e^{-|y-z|^2/(2h_t^2)}\,p_0(\dd z)}{\int e^{-|y-z'|^2/(2h_t^2)}\,p_0(\dd z')},\qquad y=e^tx,
\end{equation}
where the second expression uses $|x-e^{-t}z|^2/\sigma_t^2=|e^tx-z|^2/h_t^2$. We write $m_t(x)$ and $\Cov_t(x)$ for the mean and covariance of $\Pi_t(\cdot\mid x)$, that is, the optimal denoiser and its uncertainty; in data coordinates, $m^h(y):=\E[X_0\mid X_0+hZ=y]$, so that $m^{h_t}(e^tx)=m_t(x)$. Differentiating the formula for $p_t$ under the integral sign and dividing by $p_t$ gives $s_t(x)=\E[(e^{-t}X_0-x)/\sigma_t^2\mid X_t=x]$, and a second differentiation gives the Hessian. These are Tweedie's identities \cite{Robbins1956,Efron2011}:
\begin{equation}\label{eq:tweedie}
s_t(x)=\frac{e^{-t}m_t(x)-x}{\sigma_t^2},\qquad
\Hess\log p_t(x)=\frac{1}{\sigma_t^{2}}\Big(-\Id+\frac{\Cov_t(x)}{h_t^{2}}\Big),\qquad
\nabla m_t(x)=\frac{e^{-t}}{\sigma_t^2}\Cov_t(x).
\end{equation}
Thus the curvature of $\log p_t$ is controlled by $\Cov_t(x)/h_t^2$. Since $\supp p_0$ is compact, with diameter $D$, Popoviciu's inequality \cite{Popoviciu1935} applied along every unit direction yields $\Cov_t(x)\preceq \frac{D^2}{4}\Id$, so $s_t$ is smooth and globally Lipschitz on $\R^d$ for each $t>0$.

\emph{Reverse (generative) dynamics.} Generation runs the noising backwards. Fix a horizon $T>0$ and let $\tau\in[0,T)$ denote generation time, which corresponds to forward time $t=T-\tau$. The reverse SDE and the probability-flow ODE,
\begin{equation}\label{eq:reverse}
\dd Y_\tau=\big(Y_\tau+2s_{T-\tau}(Y_\tau)\big)\dd\tau+\sqrt2\dd W_\tau,\qquad
\dot Y_\tau=Y_\tau+s_{T-\tau}(Y_\tau),\qquad 0\le\tau<T,
\end{equation}
in which $W$ is a standard Brownian motion in $\R^d$, both have marginals $\Law(Y_\tau)=p_{T-\tau}$ when initialized at $Y_0\sim p_T$ \cite{Anderson1982,HaussmannPardoux1986,Song2021SDE}.

We label states by \emph{forward} time $t$: ``the reverse SDE started at $x$ at time $t$'' means \eqref{eq:reverse} with $Y_{T-t}=x$, whose law, indexed by forward time, does not depend on $T\ge t$. Its drift at forward time $t$ is $b_t(x):=x+2s_t(x)$, and its endpoint is $Y_{\rm end}:=\lim_{\tau\uparrow T}Y_\tau$, whose almost-sure existence is part of Theorem~\ref{thm:committor}. Since $s_t$ is Lipschitz in $x$, uniformly for $t$ in compact subintervals of $(0,\infty)$, the ODE generates a unique global diffeomorphism $\Phi_{t\to s}:\R^d\to\R^d$ from forward time $t$ to any $s\in(0,t]$. Substituting $y=e^tx$, \eqref{eq:tweedie}, and $\dd\log h_t=\dd t/\sigma_t^2$ transforms the ODE into the variance-exploding flow of \cite{Karras2022},
\begin{equation}\label{eq:VEflow}
\frac{\dd y}{\dd\log h}=y-m^h(y),
\end{equation}
integrated in the direction of decreasing $h$.

\subsection{Data models}\label{subsec:data}

We consider two classes of compactly supported data laws:
\begin{enumerate}[label=\textbf{(D\arabic*)},leftmargin=2.6em,itemsep=2pt]
\item\label{D1} \emph{Separated mixtures:} $p_0=\sum_{k=1}^Kw_k\nu_k$ with $K\ge 2$, weights $w_k>0$ summing to $1$, and probability measures $\nu_k$ having compact, pairwise disjoint supports $K_k:=\supp\nu_k$. We write $\chi\in\{1,\dots,K\}$ for the class label, so that $X_0\in K_\chi$, and $p^{(k)}_t$ for the density of $X_t$ given $\chi=k$.
\item\label{D2} \emph{Atomic data:} $p_0=\sum_{i=1}^nw_i\delta_{x_i}$ with distinct points $x_i\in\R^d$ and weights $w_i>0$, $\sum_{i=1}^n w_i=1$. This is \ref{D1} with one-point classes; equal weights give the empirical measure.
\end{enumerate}

\subsection{Committors, partitions, and observables}

\begin{definition}[Committors \cite{EVandenEijnden2006}]\label{def:committor}
For bounded measurable $g:\R^d\to\R$, the \emph{committor} is $q_g(x,t):=\E[g(Y_{\rm end})]$ for the reverse SDE started at $x$ at time $t$. Under \ref{D1} the class committors are $q_k:=q_{\one_{K_k}}$; under \ref{D2} the atom committors are $q_i:=q_{\one_{\{x_i\}}}$.
\end{definition}

\begin{definition}\label{def:interface}
Assigning each state to its most probable class gives a hard partition: under \ref{D1} (or \ref{D2}, with the atoms as classes), the \emph{Bayes cells} are $B_k(t):=\{x\in\R^d:\ q_k(x,t)>\max_{l\neq k}q_l(x,t)\}$. For $K=2$, let $\Lambda_t:=\log(q_1/q_2)$ be the committor log-odds. The \emph{interface} is $\Gamma_t:=\{\Lambda_t=0\}=\{q_1=q_2=\tfrac12\}$, and the \emph{width} of the speciation layer is
\begin{equation}\label{eq:width}
w(t):=\Prob\big(|\Lambda_t(X_t)|<\log3\big)=\Prob\big(\tfrac14<q_1(X_t,t)<\tfrac34\big),\qquad X_t\sim p_t.
\end{equation}
Thus $w(t)$ is a probability, not a geometric thickness (for which see Theorem~\ref{thm:interface}(c)).
\end{definition}

\begin{definition}[Clones \cite{Biroli2024Regimes}]\label{def:clones}
Two \emph{clones} at $(x,t)$ are two independent solutions of the reverse SDE started at $x$ at time $t$. Under \ref{D1}, the \emph{clone-agreement probability} $\phi(t):=\E\big[\sum_{k=1}^Kq_k(X_t,t)^2\big]$, $X_t\sim p_t$, is the probability that two clones started from a common state drawn from the sampler end in the same class.
\end{definition}

\begin{definition}[ODE basins and modes]\label{def:basins}
Under \ref{D1}, the ODE basin of class $k$ at time $t>0$ is 
\[
B^{\rm ODE}_k(t):=\{x\in\R^d:\ \dist(\Phi_{t\to s}(x),K_k)\to0 \, \text{ as } s\downarrow0\}.
\]
We denote by $N(t)\in\mathbb N\cup\{\infty\}$ the number of local maxima (\emph{modes}) of $p_t$ (it can be infinite, for instance for rotationally symmetric data, whose maxima need not be isolated), and set $\mathcal T_K:=\{t>0:\ N(t)=K\}$.
\end{definition}

\subsection{Problem formulation}\label{subsec:problem}

The score is exact, $s_t=\nabla\log p_t$ for the data law $p_0$ under consideration, except in Section~\ref{subsec:stability}, which quantifies the effect of an approximate score. At each time $t>0$ the reverse dynamics \eqref{eq:reverse} partition $\R^d$ by trajectory endpoints: softly through the committors, and sharply through the Bayes cells of the reverse SDE and the basins of the probability-flow ODE.

\begin{problem} \label{pb:main}
Characterize this partition, at the scale of classes, as the noise level $h_t$ decreases:
\begin{enumerate}[label=\textbf{(\alph*)},ref=(\alph*),leftmargin=2.2em,itemsep=3pt]
\item\label{pb:partition} \emph{Partitions, interfaces, and sampler comparison.} Under \ref{D1}, determine the committors $q_k(\cdot,t)$, a priori governed only by the backward Kolmogorov equation $\partial_t q=\lap q+b_t\cdot\nabla q$ with singular drift $b_t$ and datum $\one_{K_k}$ on $\supp p_0$ at $t=0$. Characterize the Bayes cells $B_k(t)$ and their interfaces (regularity, limit as $t\downarrow0$, layer width $w(t)$, and stability under score error), and compare $B_k(t)$ with the ODE basins $B^{\rm ODE}_k(t)$ and with semi-discrete optimal-transport cells.
\item\label{pb:speciation} \emph{Modes and the speciation scale.} Under \ref{D1}, describe the set $\mathcal T_K$: find an explicit noise window in which $p_t$ has exactly $K$ nondegenerate modes, decide whether $\mathcal T_K$ is an interval, and quantify how the speciation scale depends on cluster separation, diameters, and weights.
\end{enumerate}
\end{problem}

\section{Exact structure: committors are posteriors}\label{sec:exact}

Throughout this section, $p_0$ is a compactly supported probability measure on $\R^d$ with $\supp p_0\subset \bar B_{R_0}$; $p_t$ and $p^h$ denote both the laws of $X_t$ and $e^tX_t$ and their densities.
The following theorem shows that, under the exact score, the reverse SDE started at $(x,t)$ terminates almost surely with the Bayesian posterior law $\Pi_t(\cdot\mid x)$, so every committor is an explicit posterior expectation.

\begin{theorem}\label{thm:committor}
Let $0<t\le T$, $x\in\R^d$, and let $(Y_\tau)_{\tau\in[T-t,T)}$ solve the reverse SDE~\eqref{eq:reverse} with $Y_{T-t}=x$. Then:
\begin{enumerate}[label=(\roman*),leftmargin=2em]
\item $(Y_{T-s})_{s\in(0,t]}$ has the law of the mixture over $X_0\sim\Pi_t(\cdot\mid x)$ of OU bridges from $X_0$ to $x$ on $[0,t]$ (the canonical conditional law of $(X_s)_{s\in(0,t]}$ given $X_t=x$);
\item $Y_{\rm end}=\lim_{\tau\uparrow T}Y_\tau$ exists almost surely with $\Law(Y_{\rm end})=\Pi_t(\cdot\mid x)$;
\item for every bounded measurable $g$, $q_g(x,t)=\int g\,\dd\Pi_t(\cdot\mid x)=\E[g(X_0)\mid X_t=x]$. In particular, $q_g\in C^\infty(\R^d\times(0,\infty))$ solves the backward Kolmogorov equation (in forward time)
\begin{equation}\label{eq:committorPDE}
\partial_t q=\lap q+b_t\cdot\nabla q\ \ \text{on }\R^d\times(0,\infty),\quad\text{i.e.}\quad \partial_\tau\tilde q+b_{T-\tau}\cdot\nabla\tilde q+\lap\tilde q=0\ \ \text{for }\tilde q(x,\tau):=q_g(x,T-\tau).
\end{equation}
\end{enumerate}
\end{theorem}
\begin{proof}
The statement combines the time-reversal theorem \cite{Anderson1982,HaussmannPardoux1986} with Doob's $h$-transform \cite{Doob1957}; the following argument covers singular $p_0$ and every starting point. Let $k_r(y,\cdot)$ be the density of $N(e^{-r}y,\sigma_r^2\Id)$, the transition density of \eqref{eq:OU}. Fix $0<s<t$ and $\varphi\in C_b(\R^d)$, and for $t'>s$ let $v(z,t'):=\E[\varphi(X_s)\mid X_{t'}=z]=u(z,t')/p_{t'}(z)$ with $u(z,t'):=\int\varphi(y)\,p_s(y)\,k_{t'-s}(y,z)\dd y$. Both $u$ and $p$ solve the Fokker--Planck equation $\partial_{t'}f=\lap f+\nabla\!\cdot(zf)$. Hence $v$ is smooth and solves \eqref{eq:committorPDE} on $\R^d\times(s,T]$; moreover $|v|\le\sup|\varphi|$, and $v(\cdot,t')\to\varphi$ locally uniformly as $t'\downarrow s$, because $p_s$ is continuous and positive. By \eqref{eq:tweedie} and $\Cov_{t'}\preceq R_0^2\Id$, the drift $b_{T-\tau}$ is Lipschitz in $x$, uniformly for $\tau$ in compact subsets of $[T-t,T)$, so the reverse SDE has a unique strong solution on $[T-t,T)$, and by Itô's formula $v(Y_\tau,T-\tau)$, $\tau\in[T-t,T-s)$, is a bounded martingale. Letting $\tau\uparrow T-s$ gives $\E[\varphi(Y_{T-s})]=v(x,t)=\E[\varphi(X_s)\mid X_t=x]$. The same argument, started from an arbitrary state at an arbitrary time in $(s,t]$, shows that the Markov process $(Y_{T-s})_{s\in(0,t]}$ has the backward transition kernels of the forward process. Its finite-dimensional laws are therefore those of the mixture of OU bridges, and since both processes have continuous paths, their laws on $C((0,t];\R^d)$ coincide; this is (i). The bridges are continuous on $[0,t]$ with value $X_0$ at time $0$, which gives (ii) and the formula in (iii). Finally $q_g=\int g(z)\,k_t(z,\cdot)\,p_0(\dd z)/p_t$ is again a quotient of two solutions of the Fokker--Planck equation, smooth for $t>0$ by differentiation under the integral, which gives \eqref{eq:committorPDE} for $q_g$.
\end{proof}

\begin{corollary}[Explicit committors]\label{cor:explicit}
\begin{enumerate}[label=(\alph*),leftmargin=2em]
\item Under \ref{D1}, $q_k(x,t)=w_kp^{(k)}_t(x)/p_t(x)$. The pairwise log-odds $\Lambda_{kl}:=\log(q_k/q_l)$ are real-analytic on $\R^d\times(0,\infty)$, with
\begin{equation}\label{eq:gradLambda}
\nabla\Lambda_{kl}(x,t)=\frac{e^{-t}}{\sigma_t^2}\big(m^{(k)}_t(x)-m^{(l)}_t(x)\big),
\end{equation}
where $m^{(k)}_t(x)=\E[X_0\mid X_t=x,\chi=k]\in\conv K_k$ is the class-conditional posterior mean.
\item Under \ref{D2}, $q_i(x,t)=\exp\!\big(\log w_i-\frac{|e^tx-x_i|^2}{2h_t^2}\big)\big/\sum_j\exp\!\big(\log w_j-\frac{|e^tx-x_j|^2}{2h_t^2}\big)$. Hence $e^tB_i(t)$ is the open cell of $x_i$ in the power diagram of $\{x_j\}$ with weights $\psi_j=2h_t^2\log w_j$ \cite{Aurenhammer1987}. For equal weights $e^tB_i(t)$ is the Voronoi cell of $x_i$ for every $t>0$.
\item For every bounded measurable $g$, the process $\tau\mapsto q_g(Y_\tau,T-\tau)$, $\tau\in[T-t,T)$, is a bounded martingale along the reverse SDE, closed by $g(Y_{\rm end})$. Under \ref{D1}, $\int q_k(x,t)\,p_t(x)\dd x=w_k$.
\item \emph{(Persistence of class confidence.)} Under \ref{D1}, let $0<\varepsilon<\eta<1$ and $q_k(x,t)\ge1-\varepsilon$. Then the reverse SDE started at $x$ at time $t$ satisfies $\Prob\big(\inf_{T-t\le\tau<T}q_k(Y_\tau,T-\tau)\le1-\eta\big)\le\varepsilon/\eta$ and $\Prob(Y_{\rm end}\notin K_k)\le\varepsilon$. In particular, with probability at least $1-2\varepsilon$, $q_k(Y_\tau,T-\tau)>\frac12$ for all $\tau$, so that class $k$ remains the most probable class along the whole trajectory.
\end{enumerate}
\end{corollary}

\begin{proof}
(a) The supports are disjoint, so $\{X_0\in K_k\}=\{\chi=k\}$ almost surely and Bayes' rule gives the formula for $q_k$; each $p^{(k)}_t$ is positive and real-analytic in $(x,t)$, and \eqref{eq:gradLambda} is Tweedie's formula \eqref{eq:tweedie} applied to $p^{(k)}_t$. (b) is \eqref{eq:posterior} for atomic $p_0$: $q_i>q_j$ if and only if $|e^tx-x_i|^2-\psi_i<|e^tx-x_j|^2-\psi_j$. (c) By the Markov property and Theorem~\ref{thm:committor}(iii), $q_g(Y_\tau,T-\tau)=\E[g(Y_{\rm end})\mid\mathcal F_\tau]$, with $(\mathcal F_\tau)$ the filtration generated by $Y$; and $\int q_kp_t=\Prob(X_0\in K_k)=w_k$. (d) By (c), $1-q_k(Y_\tau,T-\tau)$ is a nonnegative martingale with initial value at most $\varepsilon$, so Doob's maximal inequality gives the first bound; the second is $1-q_k(x,t)\le\varepsilon$, and the last statement is the case $\eta=\frac12$ (it is void if $\varepsilon\ge\frac12$).
\end{proof}

By (c) and Jensen's inequality, $\E F\big((q_k(X_t,t))_k\big)$ with $X_t\sim p_t$ is non-increasing in $t$ for every continuous convex function $F$ on the simplex; in particular, the clone-agreement probability $\phi(t)$ of Definition~\ref{def:clones} is non-increasing.

\begin{proposition}[ODE basins carry prescribed masses]\label{prop:odemass}
Under \ref{D1}, for every $t>0$ the basins $B^{\rm ODE}_k(t)$ are Borel and pairwise disjoint, their union has full $p_t$-measure, and $p_t(B^{\rm ODE}_k(t))=w_k$.
\end{proposition}

\begin{proof}
The basins are Borel, because by continuity of the flow the convergence in Definition~\ref{def:basins} can be tested along rational $s$, and they are pairwise disjoint, because the sets $K_k$ have positive mutual distances. In data coordinates the flow \eqref{eq:VEflow} transports $p^{h_t}$ to $p^h$; since $m^h(y)\in\bar B_{R_0}$, its trajectories $y(h)$ stay in the ball $\bar B_{R(y_0)}$, $R(y_0)=\max(|y_0|,R_0)$, and move with speed at most $2R(y_0)$ in $\log h$. Let $G_\varepsilon:=\{\dist(\cdot,\supp p_0)\ge\varepsilon\}$. Since $p^h(G_\varepsilon)\le\Prob(h|Z|\ge\varepsilon)$, Tonelli's theorem and invariance give $\int p^{h_t}(\dd y_0)\int_0^{h_t}\one\{y(h)\in G_\varepsilon\}\,\dd h/h=\int_0^{h_t}p^h(G_\varepsilon)\,\dd h/h<\infty$, so almost every trajectory spends finite $\log h$-time in $G_{\varepsilon/2}$, for every $\varepsilon$ in a sequence decreasing to $0$. Such a trajectory cannot remain in $G_{\varepsilon/2}$ as $h\downarrow0$, and each passage from $\{\dist\le\varepsilon/2\}$ to $G_\varepsilon$ takes $\log h$-time at least $\varepsilon/(4R(y_0))$ inside $G_{\varepsilon/2}$; hence only finitely many passages occur, and $\dist(y(h),\supp p_0)\to0$. For small $\varepsilon$ the $\varepsilon$-neighbourhoods $U_k$ of the $K_k$ are disjoint, so the trajectory converges to a single $K_k$, as does $\Phi_{t\to s}(x)=e^{-s}y(h_s)$ for $y(h_t)=e^tx$, and $p^{h_t}(e^tB^{\rm ODE}_k(t))=\lim_{h\downarrow0}p^h(U_k)=w_k$ by dominated convergence and invariance, as $p^h\to p_0$ weakly and $p_0(\partial U_k)=0$.
\end{proof}
For three equally weighted non-collinear atoms, the following theorem shows that the probability-flow ODE yields a quadratic-cost optimal-transport partition if and only if the atoms form an equilateral triangle.

\begin{theorem}\label{thm:threeatoms}
Let $d\ge2$, $t>0$, and $p_0=\frac13(\delta_{x_1}+\delta_{x_2}+\delta_{x_3})$ with $x_1,x_2,x_3\in\R^d$ non-collinear.
\begin{enumerate}[label=(\alph*),nosep]
\item If $x_1x_2x_3$ is equilateral, then each $e^tB^{\rm ODE}_k(t)$ contains the open Voronoi cell of $x_k$, and the ODE partition is the quadratic-cost optimal-transport partition from $p^{h_t}$ to $p_0$.
\item Otherwise, the sets $e^tB^{\rm ODE}_k(t)$ do not coincide (up to Lebesgue-null sets) with the cells $L_k=\{y:|y-x_k|^2-\psi_k\le|y-x_j|^2-\psi_j\text{ for all }j\}$ of any power diagram ($\psi\in\R^3$); in particular, they do not form a quadratic-cost optimal-transport partition from $p^{h_t}$ to $p_0$. If $[x_2,x_3]$ is a longest side, $e^tB^{\rm ODE}_1(t)$ contains the open Voronoi cell of $x_1$ and a ball outside it.
\end{enumerate}
\end{theorem}

\begin{proof}
Write $h_0:=h_t$ and parametrize \eqref{eq:VEflow} by $s:=\log(h_0/h)\ge0$, so that $\dot y=m^h(y)-y$. Let $q_k$ be the posterior weights of the atoms at $(y,h)$, $\ell_{ij}:=|x_i-x_j|$ and $b_{ij}(y):=(x_i-x_j)\cdot\big(y-\frac12(x_i+x_j)\big)$, so that $\zeta_{ij}:=\log(q_i/q_j)=b_{ij}(y)/h^2$ by Corollary~\ref{cor:explicit}(b). Since $\dot b_{ij}=(x_i-x_j)\cdot(m^h-y)=\sum_lq_l\,b_{ij}(x_l)-b_{ij}$, with $b_{ij}(x_i)=-b_{ij}(x_j)=\frac12\ell_{ij}^2$ and $b_{ij}(x_k)=\frac12(\ell_{jk}^2-\ell_{ik}^2)$, we get for $\{i,j,k\}=\{1,2,3\}$
\begin{equation}\label{eq:logodds3}
\dot\zeta_{ij}=\zeta_{ij}+\frac{1}{2h^2}\Big[\ell_{ij}^2\,(q_i-q_j)+(\ell_{jk}^2-\ell_{ik}^2)\,q_k\Big],
\end{equation}
where $q_i-q_j$ has the sign of $\zeta_{ij}$. If $\zeta_{ij}$ and $\zeta_{ik}$ tend to $+\infty$ along a trajectory, then $q_i\to1$, $m^h(y)\to x_i$, and therefore $y\to x_i$.

(a) All $\ell_{ij}$ are equal, so \eqref{eq:logodds3} gives $\dot\zeta_{ij}\ge\zeta_{ij}$ wherever $\zeta_{ij}>0$. On the open Voronoi cell of $x_i$ both $\zeta_{ij}$ and $\zeta_{ik}$ are positive, hence increase to $+\infty$, and the cell lies in the basin of $x_i$. The three cells have equal $p^{h_t}$-mass by symmetry, and a power diagram whose cells have the masses of the atoms is the optimal partition \cite{AurenhammerHoffmannAronov1998}.

(b) Label the atoms so that $[x_2,x_3]$ is a longest side and $\ell_{13}<\ell_{23}$; this is possible because the triangle is not equilateral. \emph{Step 1.} On the open Voronoi cell $V$ of $x_1$ we have $q_1>q_2,q_3$, and \eqref{eq:logodds3} with $\ell_{23}\ge\ell_{12},\ell_{13}$ gives $\dot\zeta_{12}\ge\zeta_{12}>0$ and $\dot\zeta_{13}\ge\zeta_{13}>0$. Hence $V$ is forward invariant and lies in the basin of $x_1$. \emph{Step 2.} Let $b_{12}(P)=0<b_{13}(P)$. At $P$ we have $q_1=q_2$ and $\dot\zeta_{12}=(\ell_{23}^2-\ell_{13}^2)\,q_3/(2h_0^2)>0$, so the trajectory through $P$ enters $V$. By continuous dependence on the initial point, a neighbourhood of $P$ is carried into $V$; it lies in the basin of $x_1$ and contains a ball $U\subset\{b_{12}<-\eta\}$ for some $\eta>0$. \emph{Step 3.} Let $z:=\zeta_{21}$, $\zeta:=\zeta_{23}$, $\kappa:=\ell_{12}^2/(4\ell_{23}^2)\le\frac14$ and, for $h>0$, $G_h:=\{y:\ z>0,\ \zeta>0,\ e^{-\zeta}<\kappa\min(z,1)\}$, with $z$ and $\zeta$ evaluated at $(y,h)$. For $y\in G_h$ we have $q_3\le e^{-\zeta}<\kappa\min(z,1)$ and, as $\tanh\frac z2\ge\frac13\min(z,1)$, $q_2-q_1=(1-q_3)\tanh\frac z2\ge\frac14\min(z,1)$. With \eqref{eq:logodds3} this yields
\[
\dot z\ge z+\frac{\ell_{12}^2(q_2-q_1)-\ell_{23}^2\,q_3}{2h^2}\ge z,
\qquad
\dot\zeta\ge\zeta+\frac{\ell_{23}^2(q_2-q_3)-\ell_{12}^2\,q_1}{2h^2}\ge\zeta+\frac{\ell_{23}^2(q_2-q_1-q_3)}{2h^2}\ge\zeta .
\]
Thus $z$, $\zeta$ and $\zeta+\log\min(z,1)$ increase as long as $y(h)\in G_h$; hence a trajectory with $y(h_0)\in G_{h_0}$ stays in the sets $G_h$ for all $h\le h_0$, $z,\zeta\to+\infty$ along it, and $G_{h_0}$ lies in the basin of $x_2$. The open set $G_{h_0}\cap\{-\eta<b_{12}<0\}$ is nonempty: on the hyperplane $\{b_{21}=\eta/2\}$ the affine function $b_{23}$ is unbounded above, because the atoms are not collinear. \emph{Step 4.} If the basins were, up to null sets, the cells of a power diagram, then $L_1\subset\{b_{12}\ge c\}$ and $L_2\subset\{b_{12}\le c\}$ with $c:=\frac12(\psi_2-\psi_1)$. The ball $U$ gives $c<-\eta$, whereas Step~3 gives $c>-\eta$: a contradiction. Finally, optimal partitions for the quadratic cost are power diagrams \cite{AurenhammerHoffmannAronov1998}, and the last assertion follows from Steps~1 and~2.
\end{proof}

\begin{remark} 
Unlike the mass-preserving ODE basins, the reverse SDE gives a soft partition $\{q_k\}$ of masses $w_k$ (Corollary~\ref{cor:explicit}(c)) whose hard Bayes cells need not carry masses $w_k$ (Section~\ref{subsec:num-ode}); moreover, the probability-flow map is generally not optimal transport \cite{KhrulkovOseledets2022,LavenantSantambrogio2022}, agreeing with it for two atoms (Proposition~\ref{prop:odesde}) but among three equally weighted non-collinear atoms only for the equilateral triangle (Theorem~\ref{thm:threeatoms}). Note also that flow matching from a standard Gaussian source \cite{Pierret2026} is simply the flow~\eqref{eq:VEflow} started at $h=\infty$ under the reparametrization $I_u=uX_0+(1-u)Z$, $y=I_u/u$, and $h=(1-u)/u$.
\end{remark}

\section{Speciation}\label{sec:speciation}

We now analyze how density modes, clone agreement, and the SDE and ODE partitions evolve with the noise level, starting with two atoms in closed form before treating general data and $K$ separated clusters.

 \subsection{Two atoms}\label{subsec:twoatoms}

For two-atom mixtures, we can compute all speciation observables explicitly, cleanly separating the mode-splitting time, the frozen-drift threshold, and the SDE and ODE interfaces.

In the symmetric case of two equally weighted atoms, the committor, interfaces, stability thresholds, layer width, and clone agreement all take explicit closed forms.

\begin{proposition}\label{prop:twoatoms}
Let $p_0=\frac12(\delta_\mu+\delta_{-\mu})$ with $\mu\in\R^d\setminus\{0\}$, $\hat\mu:=\mu/|\mu|$, and signal-to-noise ratio $r_t:=|\mu|/h_t=e^{-t}|\mu|/\sigma_t$.
\begin{enumerate}[label=(\roman*),leftmargin=2em]
\item The SDE committor is $q_+(x,t)=\sg\big(2e^{-t}\mu\cdot x/\sigma_t^2\big)$, and for every $t>0$ both the SDE interface and the ODE basin boundary equal $\mu^\perp$.
\item $\Hess\log p_t(0)=-\sigma_t^{-2}\Id+e^{-2t}|\mu|^2\sigma_t^{-4}\,\hat\mu\otimes\hat\mu$; the density $p_t$ is bimodal if and only if $r_t>1$, i.e.\ $t<t_{\rm s}:=\frac12\log(1+|\mu|^2)$, and unimodal otherwise.
\item At $0$, the linearizations of the ODE field $x+s_t(x)$ and the frozen SDE drift $b_t$ have strictly negative eigenvalues on $\mu^\perp$, while along $\hat\mu$ they have eigenvalues $\lambda^{\rm ODE}_t=(1-\sigma_t^2)(|\mu|^2-\sigma_t^2)/\sigma_t^4$ (positive for all $t>0$ whenever $|\mu|\ge1$) and $\lambda^{\rm fr}_t=1-2\sigma_t^{-2}+2e^{-2t}|\mu|^2\sigma_t^{-4}$ (vanishing uniquely at $t_{\rm fr}:=\frac12\arcsinh(|\mu|^2)>t_{\rm s}$).
\item With $\xi\sim N(r_t,1)$, the layer width and clone-agreement probability depend on $t$ only through $r_t$:
\begin{equation}\label{eq:twoatomwidth}
w(t)=\Phi\Big(\tfrac{\log3}{2r_t}-r_t\Big)-\Phi\Big(-\tfrac{\log3}{2r_t}-r_t\Big),\qquad \phi(t)=1-2\,\E\big[\sg(2r_t\xi)\sg(-2r_t\xi)\big],
\end{equation}
and $w(t)=\frac{2}{\sqrt3}\frac{\varphi(r_t)}{r_t}\big(1+O(r_t^{-2})\big)$ as $t\downarrow0$.
\end{enumerate}
\end{proposition}

\begin{proof}
(i) The formula for $q_+$ follows from $|x-e^{-t}\mu|^2-|x+e^{-t}\mu|^2=-4e^{-t}\mu\cdot x$ and \eqref{eq:posterior}. In data coordinates $y=u\hat\mu+y_\perp$ the flow \eqref{eq:VEflow} decouples into $\frac{\dd u}{\dd(-\log h)}=|\mu|\tanh(|\mu|u/h^2)-u$ and $\frac{\dd y_\perp}{\dd(-\log h)}=-y_\perp$. So $y_\perp\to0$, and for $u>0$ the ratio $u/h$ increases as $h\downarrow0$, whence $\tanh(|\mu|u/h^2)\to1$ and $u\to|\mu|$: the basins are the open half-spaces $\{\pm\mu\cdot x>0\}$. (ii) $\Pi_t(\cdot\mid0)=p_0$, so $m_t(0)=0$ and $\Cov_t(0)=\mu\otimes\mu$ in \eqref{eq:tweedie}; moreover $p_t$ is a Gaussian in $\mu^\perp$ times a one-dimensional mixture, which is bimodal iff $r_t>1$ (Proposition~\ref{prop:RF} with $L=0$). (iii) The linearizations are $\Id+\Hess\log p_t(0)$ and $\Id+2\Hess\log p_t(0)$; use (ii) and $e^{-2t}=1-\sigma_t^2$. The equation $\lambda^{\rm fr}_t=0$ reads $\sigma_t^4-2(1+|\mu|^2)\sigma_t^2+2|\mu|^2=0$, whose root in $(0,1)$ gives $e^{2t}=|\mu|^2+\sqrt{1+|\mu|^4}>1+|\mu|^2$. (iv) With $\xi:=\hat\mu\cdot X_t/\sigma_t$ we have $\Lambda_t(X_t)=2r_t\xi$, and by symmetry we may condition on $X_0=\mu$, where $\xi\sim N(r_t,1)$. The asymptotics follow from $\int_{-a}^a\varphi(r+u)\dd u=\varphi(r)\frac{2\sinh(ra)}{r}(1+O(a^2))$ with $a=\frac{\log3}{2r}$ and $2\sinh(\frac12\log3)=\frac{2}{\sqrt3}$.
\end{proof}

\begin{figure}[t]
\centering
\includegraphics[width=0.8\textwidth]{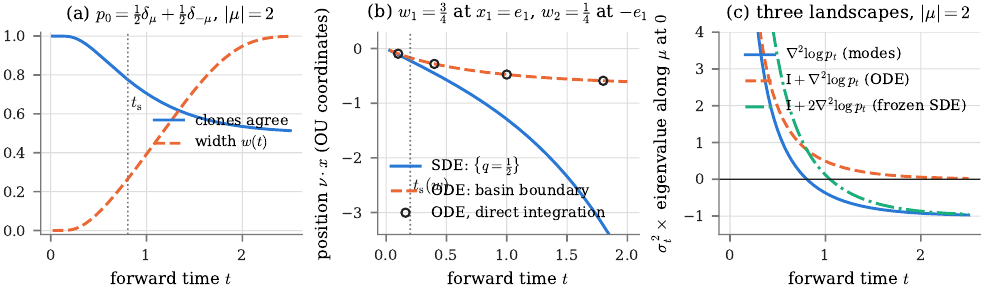}
\caption{Two atoms. (a) Equal weights, $|\mu|=2$: clone-agreement probability and width \eqref{eq:twoatomwidth} against forward time; both transition around $t_{\rm s}=\frac12\log5$. (b) Weights $\frac34$ at $e_1$ and $\frac14$ at $-e_1$: the SDE interface $e^{-t}u_S(h_t)$ and ODE basin boundary $e^{-t}u_O(h_t)$ (Proposition~\ref{prop:odesde}) in OU coordinates. Circles: thresholds found by direct ODE integration, agreeing with the quantile formula to $10^{-12}$. (c) Equal weights: eigenvalue along $\mu$ at the origin, multiplied by $\sigma_t^2$, for $\Hess\log p_t$ (modes), the ODE field $x+s_t$, and the frozen SDE drift $b_t$. Zero crossings occur at $t_{\rm s}=0.805$ (modes) and $t_{\rm fr}=\frac12\arcsinh 4=1.047$ (frozen SDE), whereas the ODE eigenvalue remains positive.}\label{fig:twoatoms}
\end{figure}

\begin{remark}[Three different ``speciation times'']\label{rem:three}
Proposition~\ref{prop:twoatoms} separates three notions that are often conflated (Figure~\ref{fig:twoatoms}c). The modes of $p_t$ split at $t_{\rm s}$. The \emph{ODE field} has a saddle at $0$ for all $t$ as soon as $|\mu|\ge1$: it undergoes no bifurcation, and its decision is encoded in the initial condition. The \emph{frozen SDE drift} acquires its saddle at $t_{\rm fr}$, and $t_{\rm fr}-t_{\rm s}\uparrow\frac12\log2$ as $|\mu|\to\infty$. For the reverse SDE the committor sharpens around $t_{\rm s}$, by \eqref{eq:twoatomwidth}. For centred data of constant norm, $0$ is a nondegenerate local maximum exactly for $t>\frac12\log(1+\lambda_1)$, with $\lambda_1$ the top eigenvalue of $\Cov(p_0)$; compare the speciation time $\frac12\log\lambda_1$ of \cite{Biroli2024Regimes}.
\end{remark}

\begin{proposition}[Two atoms, unequal weights \cite{RobertsonFryer1969,Schilling2002}]\label{prop:RF}
Let $p_0=w_1\delta_{x_1}+w_2\delta_{x_2}$ with $x_1\neq x_2$, $w_1\ge w_2>0$, $w_1+w_2=1$, $L:=\log(w_1/w_2)\ge0$, and $r_t:=|x_1-x_2|/(2h_t)$. Then $p_t$ is bimodal if and only if
\begin{equation}\label{eq:RF}
r_t>1\qquad\text{and}\qquad L<H(r_t):=2r_t\sqrt{r_t^2-1}-2\arccosh r_t,
\end{equation}
and unimodal otherwise. The function $H$ is strictly increasing on $[1,\infty)$ with $H(1)=0$ and $H(r)=2r^2-1-2\log(2r)+O(r^{-2})$ as $r\to\infty$. Consequently, $N(t)$ is non-increasing, switching from $2$ to $1$ at $t_{\rm s}(w):=\frac12\log\big(1+|x_1-x_2|^2/(4r_\star^2)\big)$ where $H(r_\star)=L$, and $h_{t_{\rm s}}^2=\frac{|x_1-x_2|^2}{2L}\big(1+O(\frac{\log L}{L})\big)$ as $w_2\downarrow0$.
\end{proposition}

\begin{proof}
Along $\nu:=(x_1-x_2)/|x_1-x_2|$, in units of $\sigma_t$ centred at the midpoint, $p_t$ is a Gaussian in $\nu^\perp$ times $f(u)=w_1\varphi(u-r)+w_2\varphi(u+r)$, $r=r_t$, and $f'$ has the sign of $g(u)+L$ with $g(u):=\log\frac{r-u}{r+u}+2ru$; the criterion follows from the variations of $g$ \cite{RobertsonFryer1969,Schilling2002}. Finally $H'(r)=4\sqrt{r^2-1}>0$ and $\arccosh r=\log(2r)-\frac{1}{4r^2}+O(r^{-4})$.
\end{proof}

\begin{proposition}[ODE versus SDE interfaces, two atoms]\label{prop:odesde}
In the setting of Proposition~\ref{prop:RF}, work in data coordinates at noise level $h>0$. Set $\nu:=(x_1-x_2)/|x_1-x_2|$, $a_i:=\nu\cdot x_i$, $D:=a_1-a_2=|x_1-x_2|$, and $c:=\frac{a_1+a_2}{2}$.
\begin{enumerate}[label=(\roman*),leftmargin=2em]
\item The Bayes cell of $x_1$ is the half-space $\{\nu\cdot y>u_S(h)\}$, where $u_S(h):=c-h^2L/D$.
\item Up to a null set, the ODE basin of $x_1$ is the half-space $\{\nu\cdot y>u_O(h)\}$, where $u_O(h)$ is the unique solution of $w_1\Phi\big(\frac{u-a_1}{h}\big)=w_2\Phi\big(\frac{a_2-u}{h}\big)$, i.e.\ the $w_2$-quantile of $\nu\cdot Y$ under $p^h$. Consequently, the ODE endpoint map coincides $p_t$-a.e.\ with the quadratic-cost optimal-transport map from $p_t$ to $p_0$.
\item If $L=0$, then $u_O\equiv u_S$. If $L>0$, then $u_S(h)<u_O(h)$ for every $h>0$, so the Bayes cell of the heavier atom strictly contains its ODE basin. Moreover, $u_O(h)-u_S(h)=\frac{hL}{2r^3}\big(1+O(r^{-2})\big)$ as $h\downarrow0$ with $r=D/(2h)$. In OU coordinates, as $t\to\infty$, $e^{-t}u_S(h_t)\to-\infty$ whereas $e^{-t}u_O(h_t)\to\Phi^{-1}(w_2)$.
\end{enumerate}
\end{proposition}

\begin{proof}
(i) is Corollary~\ref{cor:explicit}(b). (ii) The posterior weights depend on $y$ only through $u=\nu\cdot y$, so \eqref{eq:VEflow} gives a closed scalar equation for $u$, while the component of $y$ orthogonal to $\nu$ converges to that of the atoms. Scalar flows preserve order, so by Proposition~\ref{prop:odemass} the basin of $x_1$ is, up to a null set, $\{u>u^\star\}$ with $p^h(\{u>u^\star\})=w_1$. A semi-discrete optimal-transport partition for two points is a power diagram, hence a half-space orthogonal to $x_1-x_2$ of mass $w_1$: the same half-space. (iii) $F(u):=w_2\Phi(\frac{a_2-u}{h})-w_1\Phi(\frac{u-a_1}{h})$ is strictly decreasing with zero $u_O$. With $r=D/(2h)$, $A=r+\frac{L}{2r}$ and $B=r-\frac{L}{2r}$ we have $A^2-B^2=2L$, so $w_1\varphi(A)=w_2\varphi(B)$ and $F(u_S)=w_2\varphi(B)(M(B)-M(A))$, where Mills' ratio $M(z)=\Phi(-z)/\varphi(z)$ satisfies $M'=zM-1<0$ on $\R$. Hence $F(u_S)>0$, i.e.\ $u_O>u_S$, iff $L>0$. The expansion follows by linearizing $F$ at $u_S$ with $M(z)=z^{-1}-z^{-3}+O(z^{-5})$, and the limits as $t\to\infty$ from $p_t\to N(0,\Id)$ and $e^{-t}h_t^2\to\infty$.
\end{proof}

\subsection{General data: unimodality, dimension one, and failure of monotonicity}\label{subsec:modes}

By \eqref{eq:scaling}, $x$ is a critical point of $p_t$ if and only if $y=e^tx$ is a critical point of $p^{h_t}$, with $\Hess p_t(x)=e^{t(d+2)}\Hess p^{h_t}(y)$. Thus $N(t)$ is the number of local maxima of the Gaussian scale space $p^h=p_0*N(0,h^2\Id)$ of $p_0$ at scale $h=h_t$, and since $\nabla\log p^h(y)=(m^h(y)-y)/h^2$ by \eqref{eq:tweedie}, the critical points of $p^h$ are the fixed points of the mean-shift map $y\mapsto m^h(y)$ \cite{FukunagaHostetler1975,Cheng1995}. When $d=1$, $N(t)$ is non-increasing in $t$, so every $\mathcal T_K$ is an interval: $N(t)$ counts the sign changes from $+$ to $-$ of $(p^{h_t})'$, whose number is non-increasing in $h$ by the variation-diminishing property of the Gaussian kernel \cite{Karlin1968,Lindeberg1994}; for atomic data the merging times are Silverman's critical bandwidths \cite{Silverman1981}. In higher dimension this monotonicity fails (Proposition~\ref{prop:nonmonotone}).

\begin{proposition}[Sharp unimodality threshold]\label{prop:logconcave}
Let $D:=\operatorname{diam}(\supp p_0)$. If $h_t>D/2$, i.e.\ $t>t_+:=\frac12\log(1+D^2/4)$, then
\[
\Hess\log p_t(x)\preceq-\sigma_t^{-2}\Big(1-\frac{D^2}{4h_t^2}\Big)\Id\qquad\text{for all }x\in\R^d.
\]
Hence $p_t$ is strongly log-concave and has a unique, nondegenerate mode. The bound $D/2$ is sharp among data laws of diameter $D$: for two equal atoms at distance $D$, $t_{\rm s}=t_+$.
\end{proposition}

\begin{proof}
For every unit vector $v$, $v\cdot X_0$ lies in an interval of length at most $D$, so $v^\top\Cov_t(x)v\le D^2/4$ by Popoviciu's inequality \cite{Popoviciu1935}. Substituting into \eqref{eq:tweedie} yields the Hessian bound and strong log-concavity. Sharpness follows from Proposition~\ref{prop:twoatoms}(ii) with $|\mu|=D/2$.
\end{proof}

\begin{proposition}[Failure of monotonicity for $d\ge2$]\label{prop:nonmonotone}
Let $d\ge2$, and let $p_0$ be uniform on the vertices $v_1,v_2,v_3$ of an equilateral triangle of circumradius $1$ embedded in $\R^2\times\{0\}^{d-2}\subset\R^d$. Then:
\begin{enumerate}[label=(\roman*),leftmargin=2em]
\item the barycentre $0$ is a critical point of $p_t$ for all $t>0$. It is a nondegenerate local maximum for $h_t^2>\frac12$. For $h_t^2<\frac12$ it is a nondegenerate local minimum of $p_t$ restricted to the plane of the triangle: a local minimum of $p_t$ if $d=2$, a saddle if $d\ge3$;
\item at $h_t^2=\frac12$, $p_t$ possesses three additional nondegenerate local maxima along the rays $\R_+v_k$ ($k=1,2,3$). Consequently, $N(t)\ge4$ for $h_t^2\in(\frac12,\frac12+\delta)$ for some $\delta>0$, whereas $N(t)=3$ for all sufficiently small $t>0$; in particular, $N(t)$ is not non-increasing.
\end{enumerate}
\end{proposition}

\begin{proof}
(i) $\Pi_t(\cdot\mid0)=p_0$, so $m_t(0)=0$ and $\Cov_t(0)=\frac13\sum_kv_kv_k^\top=\frac12P_2$, with $P_2$ the orthogonal projection onto the plane of the triangle; by \eqref{eq:tweedie}, $\Hess\log p_t(0)=\sigma_t^{-2}(-\Id+\frac{1}{2h_t^2}P_2)$. (ii) Let $h^2=\frac12$ and $v_1=e_1$. On the ray $y=ue_1$, $u>0$, the posterior weight of $v_1$ is $w_1(u)=(1+2e^{-3u})^{-1}$, and $m^h(y)=m_1(u)e_1$ with $m_1(u)=(1-e^{-3u})/(1+2e^{-3u})$. Since $e^{3/2}>4$, we get $m_1(\frac12)>\frac12$ and $w_1(\frac12)>\frac23$; as $m_1(1)<1$, there is a critical point $u_\star e_1$ with $u_\star\in(\frac12,1)$ and $w_1:=w_1(u_\star)>\frac23$. There the posterior covariance is diagonal, with entries $\frac94w_1(1-w_1)<\frac12$ and $\frac34(1-w_1)<\frac14$ in the plane and $0$ elsewhere, so $\Hess\log p^h\prec0$. By symmetry and the implicit function theorem, three such maxima exist for $h^2$ in a neighbourhood of $\frac12$, and $N(t)=3$ for small $t$ by Theorem~\ref{thm:Kmodes} with $r=0$.
\end{proof}

\begin{figure}[t]
\centering
\includegraphics[width=0.68\textwidth]{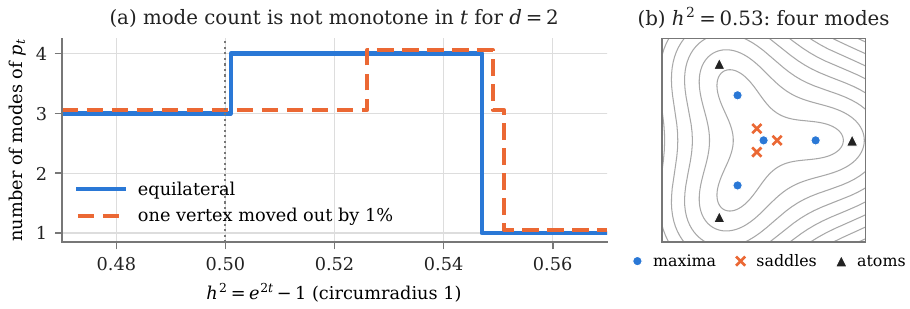}
\caption{Number of modes of $p_t$ for three equal atoms (Proposition~\ref{prop:nonmonotone}). (a) Equilateral triangle (solid) and the same triangle with one vertex moved outward by $1\%$ (dashed, shifted up by $0.06$ for visibility). (b) Level lines of $\log p^{h}$ at $h^2=0.53$, showing four maxima and three saddles.}\label{fig:modes}
\end{figure}

 Figure~\ref{fig:modes} shows the mode counts found by a multistart mean-shift and damped-Newton search. For the equilateral triangle, the ray critical points $uv_k$ ($0<u<1$) satisfy $h^2=\frac{3u}{2}\big/\log\frac{1+2u}{1-u}$, which peaks at $h^2\approx 0.54632$, where each vertex mode merges with a saddle; numerically, $N(t)$ jumps from $3$ to $4$ at $h_t^2=\frac12$ and drops to $1$ above $0.54632$. Moving one vertex outward by $1\%$ yields the sequence $3,4,3,1$ as $h_t^2$ increases, indicating that $\mathcal T_3$ need not be an interval.

\begin{remark} 
Duistermaat showed that isotropic homoscedastic Gaussian mixtures can have more modes than components for the equilateral triangle \cite{CarreiraPerpinanTR2003}, refuting the conjecture of \cite{CarreiraPerpinan2003}.
\end{remark}

\subsection{Separated mixtures: exactly \texorpdfstring{$K$}{K} modes}\label{subsec:Kmodes}

\begin{theorem}[$K$ separated clusters]\label{thm:Kmodes}
Assume \ref{D1} with $K_k\subset\bar B(c_k,r)$ ($r\ge0$), $|c_k-c_l|\ge\Delta>0$ for $k\neq l$, $D\ge\max(\Delta,\operatorname{diam}\supp p_0)$, and $\min_kw_k\ge w_{\min}>0$. Let $h:=h_t$ and $q_k:=q_k(\cdot,t)$.
\begin{enumerate}[label=(\alph*),leftmargin=2em]
\item \emph{(A dominant class at every mode.)} If $\Delta>2r$, then at every local maximum $x$ of $p_t$,
\[
\min\big(q_k(x),q_l(x)\big)\le\frac{2h^2}{(\Delta-2r)^2}\qquad\text{for all }k\neq l.
\]
\item \emph{(Exactly $K$ modes.)} Set $\varepsilon:=2(K-1)h^2/(\Delta-2r)^2$ and $\delta':=(K-1)w_{\min}^{-1}e^{-\Delta^2/(4h^2)}$. Suppose
\begin{equation}\label{eq:Kcond}
2r\le h,\qquad 4r+2\varepsilon D\le\Delta/2,\qquad \delta'<\min\Big(\varepsilon,\frac{h^2}{4D^2}\Big).
\end{equation}
Then $p_t$ has exactly $K$ local maxima, one in each ball $\mathcal B_k:=\bar B\big(e^{-t}c_k,e^{-t}(r+\varepsilon D)\big)$. Each maximum is nondegenerate, and $\log p_t$ is strictly concave on each $\mathcal B_k$.
\end{enumerate}
\end{theorem}

\begin{proof}
(a) At a local maximum $x$ of $p_t$, $\Hess\log p_t(x)\preceq0$, so $\Cov_t(x)\preceq h^2\Id$ by \eqref{eq:tweedie}. By the law of total variance, $\Cov_t=\sum_jq_jC_j+\sum_jq_j(m_j-\bar m)(m_j-\bar m)^\top$, where $C_j,m_j$ are the class-conditional posterior covariances and means and $\bar m=\sum_jq_jm_j$. Since $m_j\in\conv K_j\subset\bar B(c_j,r)$, we have $|m_k-m_l|\ge\Delta-2r>0$. Testing $\Cov_t\preceq h^2\Id$ with $v:=(m_k-m_l)/|m_k-m_l|$ and minimizing over $v\cdot\bar m\in\R$ gives
\[
h^2\ge q_k\big(v\cdot(m_k-\bar m)\big)^2+q_l\big(v\cdot(m_l-\bar m)\big)^2\ge\frac{q_kq_l}{q_k+q_l}|m_k-m_l|^2\ge\tfrac12\min(q_k,q_l)(\Delta-2r)^2.
\]

(b) Note that \eqref{eq:Kcond} implies $4r\le\Delta/2$ (hence $\Delta>2r$) and $\varepsilon\le\Delta/(4D)\le\frac14$.

\emph{Localization.} Set $\varepsilon_0:=2h^2/(\Delta-2r)^2=\varepsilon/(K-1)$. By~(a), at any local maximum $x$ at most one class has $q_j(x)>\varepsilon_0$. Since $K\ge2$ implies $K\varepsilon_0\le 2\varepsilon\le\frac12<1=\sum_jq_j$, exactly one class $k$ satisfies $q_k(x)>\varepsilon_0$, and $\sum_{j\neq k}q_j(x)\le\varepsilon$. Because $x=e^{-t}\bar m(x)$ by \eqref{eq:tweedie} and $|\bar m-m_k|\le\sum_{j\neq k}q_j|m_j-m_k|\le\varepsilon D$, we obtain $|e^tx-c_k|\le r+\varepsilon D$, i.e.\ $x\in\mathcal B_k$. The balls $\mathcal B_k$ are pairwise disjoint since $2(r+\varepsilon D)\le4r+2\varepsilon D\le\Delta/2<\Delta$.

\emph{Dominance of class $k$ on $\mathcal B_k$.} For $x\in\mathcal B_k$ and $y=e^tx$, we have $|y-z|\le2r+\varepsilon D$ for $z\in K_k$ and $|y-z|\ge\Delta-2r-\varepsilon D>0$ for $z\in K_j$ ($j\neq k$). Therefore,
\[
\frac{q_j(x)}{q_k(x)}\le\frac{w_j}{w_k}\exp\!\Big(-\frac{(\Delta-2r-\varepsilon D)^2-(2r+\varepsilon D)^2}{2h^2}\Big)\le\frac{1}{w_{\min}}e^{-\Delta^2/(4h^2)},
\]
using $(\Delta-2r-\varepsilon D)^2-(2r+\varepsilon D)^2=\Delta(\Delta-4r-2\varepsilon D)\ge\Delta^2/2$. Summing over $j\neq k$ yields $\sum_{j\neq k}q_j\le\delta'$ on $\mathcal B_k$.

\emph{Strict concavity on $\mathcal B_k$.} On $\mathcal B_k$, we have $C_k\preceq r^2\Id$, $\sum_{j\neq k}q_jC_j\preceq\delta'\frac{D^2}{4}\Id$, and $\sum_jq_j(m_j-\bar m)(m_j-\bar m)^\top\preceq\sum_{j\neq k}q_j(m_j-m_k)(m_j-m_k)^\top\preceq\delta'D^2\Id$. Using $2r\le h$ and $\delta'<h^2/(4D^2)$ gives $\Cov_t\preceq(r^2+2\delta'D^2)\Id\prec h^2\Id$, so $\Hess\log p_t\prec0$ on $\mathcal B_k$ by \eqref{eq:tweedie}.

\emph{Existence and uniqueness.} For $x\in\partial\mathcal B_k$, write $y=e^tx$ and $\varrho:=r+\varepsilon D>0$. Then $|\bar m-c_k|\le r+\delta'D$ on $\mathcal B_k$, so
\[
\nabla\log p_t(x)\cdot(y-c_k)=\frac{e^{-t}}{\sigma_t^2}\big[(\bar m-c_k)\cdot(y-c_k)-\varrho^2\big]\le\frac{e^{-t}\varrho}{\sigma_t^2}(r+\delta'D-\varrho)<0,
\]
because $\delta'<\varepsilon$. Thus $\nabla\log p_t$ points strictly inward on $\partial\mathcal B_k$, so $p_t$ attains its maximum over $\mathcal B_k$ at an interior critical point, which is unique and nondegenerate by strict concavity. Combined with localization, $p_t$ has exactly $K$ local maxima.
\end{proof}

\begin{corollary}\label{cor:window}
Under the assumptions of Theorem~\ref{thm:Kmodes} with $r\le\Delta/16$, set
\[
h_+:=\min\!\Bigg(\frac{\Delta^{3/2}}{\sqrt{21(K-1)D}},\;\frac{\Delta}{\sqrt{8\log\!\big(16(K-1)D^2/(w_{\min}\Delta^2)\big)}}\Bigg).
\]
Whenever $2r\le h_+$, the density $p_t$ has exactly $K$ nondegenerate modes for all $h_t\in[2r,h_+]$ and a single mode for $h_t>D/2$ (Proposition~\ref{prop:logconcave}), so the passage from $K$ modes to one takes place in $(h_+,D/2]$.
\end{corollary}

\begin{proof}
Fix $h_t\in[2r,h_+]$. Since $r\le\Delta/16$, we have $\Delta-2r\ge\frac78\Delta$ and $\varepsilon\le\frac{128}{49}(K-1)h_t^2/\Delta^2$. By $h_t\le h_+$, this yields $2\varepsilon D\le\frac{1024}{49\cdot21}\frac{\Delta}{4}<\frac{\Delta}{4}$, so $4r+2\varepsilon D\le\Delta/2$. Because $D\ge\Delta$ and $K\ge2$, we have $\varepsilon\ge2h_t^2/\Delta^2>h_t^2/(4D^2)$, so the last condition in \eqref{eq:Kcond} reduces to $\delta'<h_t^2/(4D^2)$, i.e.\ $u>\log(Mu)$ where $u:=\Delta^2/(4h_t^2)$ and $M:=16(K-1)D^2/(w_{\min}\Delta^2)\ge32$. Since $u\mapsto u-\log u$ increases on $[1,\infty)$ and $M>2\log M$ for all $M>0$, this holds whenever $u\ge2\log M$, i.e.\ $h_t\le\Delta/\sqrt{8\log M}$, which is ensured by $h_t\le h_+$. Theorem~\ref{thm:Kmodes}(b) then gives $K$ nondegenerate modes.
\end{proof}

When $K$ and $D/\Delta$ are bounded, $h_+$ is of order $\Delta/\sqrt{\log(1/w_{\min})}$, matching the sharp two-atom scaling of Proposition~\ref{prop:RF} (for two equal atoms at distance $D$, $h_+=D/\sqrt{8\log32}\approx0.190\,D$, certifying $38\%$ of the exact bimodal range $h_t<D/2$).
\begin{proposition}[Speciation window for clone agreement]\label{prop:clonewindow}
Under the assumptions of Theorem~\ref{thm:Kmodes} with $\Delta>2r$, for every $t>0$ and $h=h_t$,
\begin{equation}\label{eq:clonewindow}
\Big(1-\sum_{k=1}^Kw_k^2\Big)\,e^{-D^2/(4h^2)}\ \le\ 1-\phi(t)\ \le\ (K-1)\,e^{-(\Delta-2r)^2/(8h^2)}.
\end{equation}
Consequently, if $\phi(t)=\theta\in(\sum_kw_k^2,1)$, then
$\ \frac{(\Delta-2r)^2}{8\log((K-1)/(1-\theta))}\le h_t^2\le\frac{D^2}{4\log((1-\sum_kw_k^2)/(1-\theta))}$.
\end{proposition}

\begin{proof}
Let $a_k:=w_kp^{(k)}_t$, so that $q_k=a_k/p_t$ and $1-\phi(t)=\sum_{k\neq l}\int a_ka_l/p_t\dd x$. \emph{Lower bound.} By the Cauchy--Schwarz inequality, $(\int\sqrt{a_ka_l}\dd x)^2\le\int(a_ka_l/p_t)\dd x\int p_t\dd x=\int(a_ka_l/p_t)\dd x$. In data coordinates the class densities are the mixtures $\int N(z,h^2\Id)\,\nu_k(\dd z)$, and two Gaussians $N(z,h^2\Id)$, $N(z',h^2\Id)$ have Bhattacharyya coefficient $e^{-|z-z'|^2/(8h^2)}\ge e^{-D^2/(8h^2)}$; since $(a,b)\mapsto\sqrt{ab}$ is concave, Jensen's inequality gives $\int(p^{(k)}_tp^{(l)}_t)^{1/2}\dd x\ge e^{-D^2/(8h^2)}$. Hence $\int a_ka_l/p_t\dd x\ge w_kw_le^{-D^2/(4h^2)}$, and $\sum_{k\neq l}w_kw_l=1-\sum_kw_k^2$. \emph{Upper bound.} Since $\sum_kq_k^2\ge\max_kq_k^2$ and $1-m^2\le2(1-m)$, $1-\phi(t)\le2\,\E[1-\max_kq_k(X_t,t)]$, twice the Bayes error for classifying $\chi$ from $X_t$, which is at most the error of any classifier. Assign $y=e^tX_t=X_0+hZ$ to the nearest centre. If $\chi=k$ and $c_l$ is at least as near as $c_k$, then $u:=(c_l-c_k)/|c_l-c_k|$ satisfies $hZ\cdot u\ge\frac12|c_l-c_k|-(X_0-c_k)\cdot u\ge\frac12(\Delta-2r)$. A union bound over $l\neq k$ gives the error bound $(K-1)(1-\Phi(\frac{\Delta-2r}{2h}))$, and $1-\Phi(x)\le\frac12e^{-x^2/2}$ for $x\ge0$. The last assertion follows by solving \eqref{eq:clonewindow} for $h$.
\end{proof}

Thus clone agreement passes from $\sum_kw_k^2$ to $1$ at noise levels between multiples of $\Delta-2r$ and of $D$, the range in which the modes of $p_t$ merge (Corollary~\ref{cor:window}). The exponent of the upper bound is sharp: for two equal atoms, \eqref{eq:twoatomwidth} gives $1-\phi(t)\sim\frac\pi2\varphi(r_t)/r_t$ as $r_t\to\infty$.

\section{Interfaces}\label{sec:interfaces}

Throughout this section we take $K=2$ under \ref{D1}, write $\Lambda_t:=\log(q_1/q_2)=\log\frac{w_1p^{(1)}_t}{w_2p^{(2)}_t}$, and denote by $m^{(k)}_t,C^{(k)}_t$ the class-conditional posterior means and covariances. We say that the two classes are \emph{convexly separated} if there exist a unit vector $v\in\R^d$ and $g>0$ such that
\begin{equation}\label{eq:sep}
v\cdot(z_1-z_2)\ge g\qquad\text{for all }z_1\in K_1,\ z_2\in K_2,
\end{equation}
and we set $D:=\max_{z_1\in K_1,\,z_2\in K_2}|z_1-z_2|\ge g$. The \emph{bisector} of the supports is $\Bis:=\{y\in\R^d:\dist(y,K_1)=\dist(y,K_2)\}$.

\subsection{Regularity, sharpness, and the small-noise limit}\label{subsec:regularity}

\begin{theorem}[Shape of the SDE interface]\label{thm:interface}
Assume \eqref{eq:sep}. For every $t>0$:
\begin{enumerate}[label=(\alph*),leftmargin=2em]
\item $v\cdot\nabla\Lambda_t\ge e^{-t}g/\sigma_t^2$ and $|\nabla\Lambda_t|\le e^{-t}D/\sigma_t^2$ on $\R^d$.
\item $\Gamma_t$ is the graph $\{z+\gamma_t(z)v:\ z\in v^\perp\}$ of a real-analytic function $\gamma_t:v^\perp\to\R$ with $\operatorname{Lip}(\gamma_t)\le\sqrt{D^2-g^2}/g$. In particular, $\Gamma_t$ is a connected real-analytic hypersurface for every $t>0$, both before and after speciation.
\item For every $a>0$ and $z\in v^\perp$, the set $\{x\in z+\R v:\ |\Lambda_t(x)|<a\}$ is an open interval containing $z+\gamma_t(z)v$, of length between $2a\,e^t\sigma_t^2/D$ and $2a\,e^t\sigma_t^2/g$. In data coordinates the length lies between $2ah_t^2/D$ and $2ah_t^2/g$.
\item In data coordinates the interface $e^t\Gamma_t$ is the graph of $\gamma^\sharp_t(z):=e^t\gamma_t(e^{-t}z)$. As $t\downarrow0$, $\gamma^\sharp_t\to\gamma_0$ locally uniformly on $v^\perp$, where $\Bis$ is the graph of $\gamma_0$ (with $\operatorname{Lip}(\gamma_0)\le\sqrt{D^2-g^2}/g$). If moreover there exist $c>0$ and $m\ge0$ such that $\nu_k(B(z,\rho))\ge c\rho^m$ for all $z\in K_k$, $\rho\in(0,1]$, and $k=1,2$, then for every $R>0$ there exists $C_R>0$ such that $\sup_{|z|\le R}|\gamma^\sharp_t(z)-\gamma_0(z)|\le C_Rh_t^2\log(1/h_t)$ for all sufficiently small $t>0$; if $m=0$, as for atomic classes, the bound improves to $C_Rh_t^2$.
\end{enumerate}
\end{theorem}

\begin{proof}
(a) By \eqref{eq:gradLambda}, $\nabla\Lambda_t=e^{-t}\sigma_t^{-2}(m^{(1)}_t-m^{(2)}_t)$, and every $u\in\conv K_1-\conv K_2$ satisfies $v\cdot u\ge g$ and $|u|\le D$.

(b) Integrating (a) along $z+\R v$ shows that $s\mapsto\Lambda_t(z+sv)$ is strictly increasing with range $\R$. The analytic implicit function theorem gives $\gamma_t$ with $\nabla\gamma_t=-P_{v^\perp}\nabla\Lambda_t/(v\cdot\nabla\Lambda_t)$, and $|P_{v^\perp}u|/(v\cdot u)\le\sqrt{D^2-g^2}/g$ for $u\in\conv K_1-\conv K_2$. (c) follows from (a) by the mean value theorem along the line.

(d) In data coordinates, $\Lambda=\log\frac{w_1}{w_2}+\log I_1-\log I_2$ with $I_k(y):=\int e^{-|y-z|^2/(2h^2)}\nu_k(\dd z)$. Let $d_k:=\dist(\cdot,K_k)$ and let $\pi_k(y)\in K_k$ be a nearest point. Restricting the integral to $B(\pi_k(y),\rho)$ gives
\[
0\le-2h^2\log I_k(y)-d_k(y)^2\le2\rho\,d_k(y)+\rho^2+2h^2\log\frac{1}{\nu_k(B(\pi_k(y),\rho))}\qquad(\rho>0),
\]
and $\inf_{z\in K_k}\nu_k(B(z,\rho))>0$ by compactness. Hence $-2h^2\Lambda=G+E_h$ with $G:=d_1^2-d_2^2$ and $E_h\to0$ locally uniformly; under the density condition, $\rho=h^2$ gives $|E_h|\le C_Rh^2\log(1/h)$ on $B_R$, and $\rho\downarrow0$ gives $|E_h|\le C_Rh^2$ if $m=0$. Next, $G$ decreases along $v$ at rate at least $2g$: for $s>0$ and $\pi_2':=\pi_2(y+sv)$, expanding $d_1(y+sv)^2\le|y+sv-\pi_1(y)|^2$ and $d_2(y)^2\le|(y+sv-\pi_2')-sv|^2$ and adding the two inequalities gives
\[
G(y+sv)-G(y)\le2s\,v\cdot(\pi_2'-\pi_1(y))\le-2gs .
\]
Therefore $\{G=0\}=\Bis$ is a graph $\gamma_0$ over $v^\perp$, and $2g\,|\gamma^\sharp_t(z)-\gamma_0(z)|\le|E_h(z+\gamma^\sharp_t(z)v)|$. Taking $\rho=g/2$ above gives $|E_h(y)|\le C+g\max_kd_k(y)$ for $h\le1$, which bounds $|\gamma^\sharp_t(z)|$ uniformly for $|z|\le R$ and $h_t\le1$. Both claims then follow on a fixed ball, and the Lipschitz bound passes to the limit $\gamma_0$.
\end{proof}

No restriction $t<t_{\rm s}$ is needed in Theorem~\ref{thm:interface}: what changes at speciation is not the regularity of $\Gamma_t$ but the \emph{sharpness} of the committor across it, governed by $|\nabla\Lambda_t|$.

\subsection{Width of the speciation layer}

\begin{theorem}[Width]\label{thm:width}
Let $f_t$ and $f^{(k)}_t$ denote the marginal and class-conditional densities of $\Lambda_t(X_t)$. Because $K_1\cap K_2=\emptyset$, $\Lambda_t$ is real-analytic and nonconstant, so $\nabla\Lambda_t\neq0$ almost everywhere and the density identities in~(a)--(b) hold for almost every $\lambda\in\R$. For~(c) and \eqref{eq:hessbound}, assume also \eqref{eq:sep}, under which the coarea integral in~(a) defines a continuous, strictly positive density $f_t(\lambda)$ for all $\lambda\in\R$ (Appendix~\ref{app:width}).
\begin{enumerate}[label=(\alph*),leftmargin=2em]
\item \emph{(Coarea.)} $f_t(\lambda)=\int_{\{\Lambda_t=\lambda\}}\frac{p_t}{|\nabla\Lambda_t|}\dd\mathcal H^{d-1}$ and $w(t)=\int_{-\log3}^{\log3}f_t(\lambda)\dd\lambda$.
\item \emph{(Calibration.)} $w_1f^{(1)}_t(\lambda)=e^{\lambda}w_2f^{(2)}_t(\lambda)$, so $f_t(\lambda)=w_2f^{(2)}_t(\lambda)(1+e^\lambda)$.
\item \emph{(Symmetric problems.)} Suppose $w_1=w_2$ and there exists a Euclidean isometry $R$ of $\R^d$ with $R_\#\nu_1=\nu_2$ and $R_\#\nu_2=\nu_1$. Then $f_t$ is even, $f^{(2)}_t(\lambda)=f_t(0)\,e^{\beta_t(\lambda)-\lambda/2}$ for a continuous even function $\beta_t$ with $\beta_t(0)=0$, and
\begin{equation}\label{eq:widthformula}
w(t)=f_t(0)\int_{-\log3}^{\log3}e^{\beta_t(\lambda)}\cosh\tfrac\lambda2\dd\lambda,\qquad
\Big|\log\frac{w(t)}{\frac{4}{\sqrt3}f_t(0)}\Big|\le\bar\beta_t:=\max_{|\lambda|\le\log3}|\beta_t(\lambda)|.
\end{equation}
For two equal atoms, and more generally for $X_0=c+Sav+U$ with $c\in\R^d$, a unit vector $v$, $a>0$, a sign $S=\pm1$ with equal probabilities labelling the class, and $U$ independent of $S$ with any compactly supported law on $v^\perp$ (so that \eqref{eq:sep} holds with $g=2a$), we have $\beta_t(\lambda)=-\lambda^2/(8r_t^2)$ identically, with $r_t:=a/h_t$; hence $\bar\beta_t=(\log3)^2/(8r_t^2)\to0$ as $t\downarrow0$, and $w(t)$ is given by \eqref{eq:twoatomwidth}.
\item At every regular point of $\Gamma_t$, set $n:=\nabla\Lambda_t/|\nabla\Lambda_t|$. Then
\begin{equation}\label{eq:hessgamma}
n^\top\big(\Hess\log p_t+\sigma_t^{-2}\Id\big)n=\tfrac14|\nabla\Lambda_t|^2+\frac{e^{-2t}}{2\sigma_t^4}\,n^\top\big(C^{(1)}_t+C^{(2)}_t\big)n\ge\tfrac14|\nabla\Lambda_t|^2.
\end{equation}
Consequently, under the separation and symmetry hypotheses of~(c),
\begin{equation}\label{eq:hessbound}
w(t)\ge\frac{2}{\sqrt3}\,e^{-\bar\beta_t}\int_{\Gamma_t}\frac{p_t}{\sqrt{\lambda_{\max}(\Hess\log p_t)+\sigma_t^{-2}}}\dd\mathcal H^{d-1}.
\end{equation}
Equality is asymptotic for two equal atoms as $t\downarrow0$.
\end{enumerate}
\end{theorem}

\begin{proof}
(a) If $\Lambda_t$ were constant, then $p^{(1)}_t=p^{(2)}_t$, hence $\nu_1=\nu_2$ because convolution with a Gaussian is injective, contradicting $K_1\cap K_2=\emptyset$. So the real-analytic field $\nabla\Lambda_t$ vanishes only on a null set, and the coarea formula applies. (b) For any Borel set $A\subset\R$,
\[
w_1\Prob\big(\Lambda_t(X_t)\in A\mid X_0\in K_1\big)=\E\big[q_1(X_t,t)\one_A(\Lambda_t(X_t))\big]=\E\big[\sg(\Lambda_t(X_t))\one_A(\Lambda_t(X_t))\big],
\]
and similarly for class~$2$ with $\sg(-\Lambda_t)$. Thus $w_1f^{(1)}_t(\lambda)=\sg(\lambda)f_t(\lambda)$ and $w_2f^{(2)}_t(\lambda)=\sg(-\lambda)f_t(\lambda)$ a.e., giving~(b) since $\sg(\lambda)/\sg(-\lambda)=e^\lambda$. (c) In data coordinates $Y:=e^tX_t=X_0+h_tZ$, the isometry $R$ preserves the law $p^{h_t}$ and swaps the two components, so the log-odds in data coordinates, $\Lambda^\sharp(y):=\Lambda_t(e^{-t}y)$, satisfy $\Lambda^\sharp\circ R=-\Lambda^\sharp$. Hence $\Lambda_t(X_t)\overset{d}{=}-\Lambda_t(X_t)$ and $f_t$ is even. By~(b) with $w_2=\frac12$, $f^{(2)}_t(\lambda)=2\sg(-\lambda)f_t(\lambda)=f_t(\lambda)e^{-\lambda/2}/\cosh\frac\lambda2$. Defining $\beta_t(\lambda):=\log\big(f_t(\lambda)/(f_t(0)\cosh\frac\lambda2)\big)$ yields an even continuous function with $\beta_t(0)=0$; integrating $f_t$ over $[-\log3,\log3]$ using $\int_{-\log3}^{\log3}\cosh\frac\lambda2\dd\lambda=4\sinh(\frac12\log3)=4/\sqrt3$ gives \eqref{eq:widthformula}. In the example, the two class densities in data coordinates are $\varphi((u\mp a)/h)/h$, $u:=v\cdot(y-c)$, times a common function of the component of $y$ in $v^\perp$. Hence $\Lambda^\sharp=2au/h^2$, and $\Lambda_t(X_t)=2r_t\xi$ with $\xi\sim\frac12N(r_t,1)+\frac12N(-r_t,1)$, so $f_t(\lambda)=\frac{\varphi(r_t)}{2r_t}e^{-\lambda^2/(8r_t^2)}\cosh\frac\lambda2$. (d) On $\Gamma_t$ we have $q_1=q_2=\frac12$, so the law of total covariance gives $\Cov_t=\frac12(C^{(1)}_t+C^{(2)}_t)+\frac14(m^{(1)}_t-m^{(2)}_t)^{\otimes2}$. Substituting into \eqref{eq:tweedie} and \eqref{eq:gradLambda} yields \eqref{eq:hessgamma}. Hence $\frac12|\nabla\Lambda_t|\le(\lambda_{\max}(\Hess\log p_t)+\sigma_t^{-2})^{1/2}$ on $\Gamma_t$, and combining~(a) with \eqref{eq:widthformula} gives \eqref{eq:hessbound}. For two equal atoms, $C^{(1)}_t=C^{(2)}_t=0$ and $\bar\beta_t\to0$, making \eqref{eq:hessbound} asymptotically sharp as $t\downarrow0$.
\end{proof}

\begin{remark}[Comparison with a saddle-curvature heuristic]\label{rem:widthlaw}
A natural saddle-point heuristic estimates the width by $w\asymp\lambda_{\max}^{-1/2}\int_{\Gamma_t}p_t\dd\mathcal H^{d-1}$, with $\lambda_{\max}$ the unstable eigenvalue of $\Hess\log p_t$ at the saddle. By \eqref{eq:hessgamma}, the curvature governing the layer thickness is instead the \emph{between-class} part $\frac14|\nabla\Lambda_t|^2$. Near speciation $\lambda_{\max}\downarrow0$, so the heuristic diverges while $w(t)$ stays bounded: for two equal atoms it gives $\frac{2}{\sqrt3}\varphi(r_t)/\sqrt{r_t^2-1}$ instead of $\frac{2}{\sqrt3}\varphi(r_t)/r_t\,(1+O(r_t^{-2}))$. Section~\ref{subsec:num-width} shows the same discrepancy for interlocking arcs (Figure~\ref{fig:width}).
\end{remark}

\begin{remark}[What free-boundary theory contributes]\label{rem:FB}
For the exact score at positive noise, the interfaces are level sets of explicit real-analytic functions; no singular free boundary can form, and the monotonicity formulas of Alt--Caffarelli--Friedman type \cite{ACF1984}, central to two-phase free-boundary and optimal-partition problems \cite{ContiTerraciniVerzini2005,CaffarelliLin2008,HHT2009,Boz2022Singular}, are not needed. Free-boundary and variational problems do enter in three ways. (i) \emph{Zero-noise limit:} in data coordinates, with $p^{(k),h}:=\nu_k*N(0,h^2\Id)$, the function $v_k:=-2h^2\log(w_kp^{(k),h})$ converges by the Laplace principle \cite{Varadhan1967} to $\dist(\cdot,K_k)^2$, so the limiting partition is the generalized Voronoi diagram of the supports. (ii) \emph{Approximate scores:} a learned score is in general neither exact nor a gradient, so the committor solves \eqref{eq:committorPDE} with a non-explicit drift and its interfaces pose a genuine PDE problem (Sections~\ref{subsec:stability} and~\ref{subsec:num-fd}). (iii) \emph{Frozen problem:} for fixed $t$, $b_t=\nabla V_t$ with $V_t:=\frac12|x|^2+2\log p_t$, so the frozen generator is symmetric in $L^2(e^{V_t}\dd x)$ and its committors minimize $\int|\nabla q|^2e^{V_t}\dd x$, a weighted Dirichlet partition problem \cite{Boz2015Partitions,HHT2009}; the frozen committor is not that of the time-dependent dynamics (Proposition~\ref{prop:twoatoms}(ii)--(iii)).
\end{remark}

\subsection{Stability under score error}\label{subsec:stability}

Let $\hat s:\R^d\times(0,T]\to\R^d$ be an approximate score such that, uniformly for $u$ in compact subsets of $(0,T]$, $\hat s_u$ is Lipschitz in $x$ and $\hat s_u-s_u$ is bounded. For $0<t_0<t\le T$ let $\Prob_x$ and $\hat\Prob_x$ be the laws on $C([T-t,T-t_0];\R^d)$ of the reverse SDE in \eqref{eq:reverse} and of the same SDE with $\hat s$ in place of $s$ (the \emph{learned sampler}), both started at $x$ at time $t$, with expectations $\E_x$, $\hat\E_x$ and canonical process $Y$. Let $\delta(t_0,t)^2:=\int_{t_0}^{t}\int_{\R^d}|\hat s_u-s_u|^2\,p_u\dd x\dd u$ be the score-matching error integrated over the noise levels between $t_0$ and $t$, with the weighting of maximum-likelihood training \cite{SongDurkan2021}.

\begin{theorem}\label{thm:stability}
For every measurable $\Psi:C([T-t,T-t_0];\R^d)\to[0,1]$,
\[
\int_{\R^d}\big|\hat\E_x\Psi-\E_x\Psi\big|\,p_t(x)\dd x\le\frac{\delta(t_0,t)}{\sqrt2}.
\]
In particular, under \ref{D1} the \emph{hybrid committors} $\hat q_k(x,t):=\hat\E_x[q_k(Y_{T-t_0},t_0)]$ (learned sampler down to the cutoff $t_0$, exact posterior thereafter) satisfy $\E|\hat q_k(X_t,t)-q_k(X_t,t)|\le\delta(t_0,t)/\sqrt2$ for $X_t\sim p_t$.
\end{theorem}

\begin{proof}
On $[T-t,T-t_0]$ both drifts are Lipschitz and they differ by the bounded field $2(\hat s-s)$, so Girsanov's theorem applies and, the diffusion coefficient being $\sqrt2$, the relative entropy is $\mathrm{KL}(\Prob_x\,\|\,\hat\Prob_x)=\E_x\int_{T-t}^{T-t_0}|\hat s_{T-\tau}-s_{T-\tau}|^2(Y_\tau)\dd\tau$. By Pinsker's inequality, $|\hat\E_x\Psi-\E_x\Psi|\le(\frac12\mathrm{KL}(\Prob_x\,\|\,\hat\Prob_x))^{1/2}$. Integrating against $p_t$, using Jensen's inequality and $\Law(Y_\tau)=p_{T-\tau}$ for the exact sampler started from $p_t$, gives the bound. For the committors take $\Psi=q_k(Y_{T-t_0},t_0)$ and use Corollary~\ref{cor:explicit}(c).
\end{proof}

\begin{corollary}[Stability of the Bayes partition]\label{cor:stability}
Let $K=2$ under \ref{D1}, $\hat B_1(t):=\{\hat q_1(\cdot,t)>\frac12\}$ and $\delta:=\delta(t_0,t)$. Then $p_t(\hat B_1(t)\triangle B_1(t))\le w(t)+2\sqrt2\,\delta$. If moreover \eqref{eq:sep} holds and $\delta\le0.9\,F_t$ with $F_t:=\max_{|\lambda|\le\log3}f_t(\lambda)$, then $p_t(\hat B_1(t)\triangle B_1(t))\le5.5\sqrt{F_t\,\delta}$.
\end{corollary}

\begin{proof}
Let $\eta\in(0,\frac14]$. Where $\hat q_1-\frac12$ and $q_1-\frac12$ differ in sign, either $|\hat q_1-q_1|\ge\eta$, which has $p_t$-probability at most $\delta/(\sqrt2\eta)$ by Theorem~\ref{thm:stability} and Markov's inequality, or $|q_1-\frac12|<\eta$, that is, $|\Lambda_t|<\ell(\eta):=\log\frac{1+2\eta}{1-2\eta}\le\frac{16}3\eta$. For $\eta=\frac14$ the latter event has probability $w(t)$. Under \eqref{eq:sep}, $f_t$ is continuous (Appendix~\ref{app:width}), so its probability is at most $\frac{32}3\eta F_t$, and minimizing $\frac{32}3\eta F_t+\delta/(\sqrt2\eta)$ over $\eta$ gives $2(\frac{32}{3\sqrt2}F_t\delta)^{1/2}\le5.5\sqrt{F_t\delta}$; the minimizer is at most $\frac14$ when $\delta\le0.9F_t$.
\end{proof}

Thus the hybrid and the exact partition can disagree, beyond a set of mass $2\sqrt2\,\delta$, only inside the speciation layer; since $F_t\le\frac12e^{2\bar\beta_t}w(t)$ in the setting of Theorem~\ref{thm:width}(c), the disagreement has mass at most of order $\sqrt{w(t)\,\delta}$. Girsanov bounds of this type are standard for the sampling error of diffusion models \cite{ChenChewi2023}. The cutoff cannot be removed without further information. If the learned sampler is continued below $t_0$ and then ends in class $k$ with probability $r_k(z)$ from the state $z$ at time $t_0$, then its committor $\tilde q_k(x,t):=\hat\E_x[r_k(Y_{T-t_0})]$ satisfies $\E|\tilde q_k(X_t,t)-q_k(X_t,t)|\le\delta/\sqrt2+\epsilon_k$ with $\epsilon_k:=\int|r_k-q_k(\cdot,t_0)|\dd\hat\mu_{t_0}$, where $\hat\mu_{t_0}$ is the law at time $t_0$ of the learned sampler started from $p_t$, which the score error on $[t_0,t]$ does not control.

\section{Numerical experiments}\label{sec:numerics}

\subsection{Finite-difference committors and score error}\label{subsec:num-fd}

We begin with a planar test problem consisting of two interlocking semicircular arcs (``moons''), $A_1:=\{(\cos\vartheta,\sin\vartheta):\vartheta\in[0,\pi]\}-(\frac12,\frac14)$ and $A_2:=-A_1$, each discretized by $200$ equally spaced atoms with total weight $\frac12$. The arcs are separated by distance $0.5$, but their convex hulls overlap, so \eqref{eq:sep} fails and this example tests the theory beyond Theorem~\ref{thm:interface}.

\begin{figure}[t]
\centering
\includegraphics[width=0.8\textwidth]{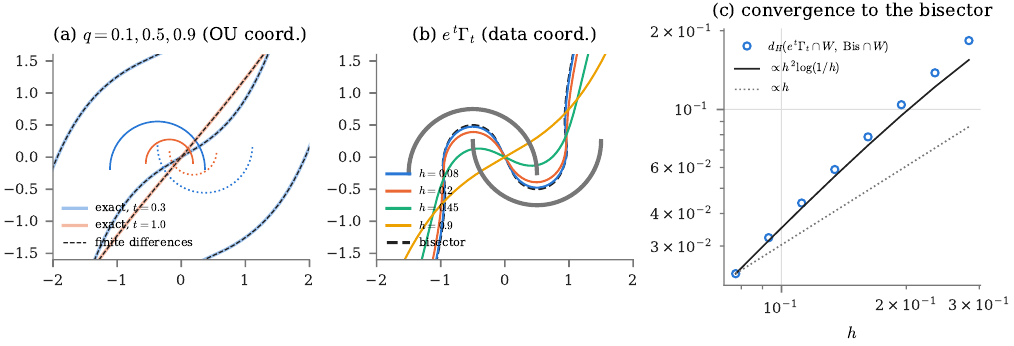}
\caption{Interlocking arcs. (a) Level sets $q_1=0.1,0.5,0.9$ of the committor at $t=0.3$ and $t=1$ (OU coordinates): exact posterior (thick, light) and finite differences (dashed); the arcs are drawn scaled by $e^{-t}$. (b) Interfaces $e^t\Gamma_t$ in data coordinates for $h=0.08,0.2,0.45,0.9$, together with the bisector of the arcs (dashed). (c) Hausdorff distance to the bisector in $W=[-2,2]\times[-1.6,1.6]$, with reference slopes; the decay matches the $h^2\log(1/h)$ rate of Theorem~\ref{thm:interface}(d) even though \eqref{eq:sep} fails.}\label{fig:moons}
\end{figure}

To verify that the exact posterior propagates according to \eqref{eq:committorPDE}, we solve the equation in forward time from the exact profile at $t_0=0.02$ using Strang splitting (exact cosine-basis heat half-steps with Neumann boundary conditions around a semi-Lagrangian midpoint advection step with cubic-spline interpolation), with drift computed from \eqref{eq:posterior}. On $[-4.5,4.5]\times[-4,4]$ with $\Delta x=0.025$ and $\Delta t=\min(0.025\,\sigma_t^2,0.02)$, the maximum error against Theorem~\ref{thm:committor} on $|x_1|<3$, $|x_2|<2.5$ is $2.4\times10^{-4}$ at $t=0.1$ and below $10^{-4}$ for $0.3\le t\le1.5$; halving the time step reduces the error by factors of $1.5$--$2$ for $t\le1$ (Figure~\ref{fig:moons}a).

To test the score-error bound of Theorem~\ref{thm:stability}, we take $d=1$, $p_0=\frac12(\delta_{-1}+\delta_1)$, and a constant perturbation $\hat s:=s+c$, so that $\delta(t_0,t)^2=c^2(t-t_0)$. We integrate \eqref{eq:committorPDE} with drift $x+2\hat s_t$ from the exact profile at $t_0=0.05$ to $t=0.6$ by Crank--Nicolson on $[-8,8]$ ($\Delta x=0.005$, $\Delta t=2.5\times10^{-4}$, boundary values $0$ and $1$). For $c=0.02,0.1,0.25$, the hybrid committor $\hat q_+$ of $+1$ satisfies $\E|\hat q_+-q_+|=0.0058,0.029,0.073$, safely below $\delta/\sqrt2=0.0105,0.052,0.131$, and the partitions disagree on sets of mass $0.0078,0.039,0.097$ (for $c=0$ the scheme reproduces $q_+$ to $1.3\times10^{-6}$).

\subsection{Interfaces and widths}\label{subsec:num-width}

\begin{figure}[t]
\centering
\includegraphics[width=0.86\textwidth]{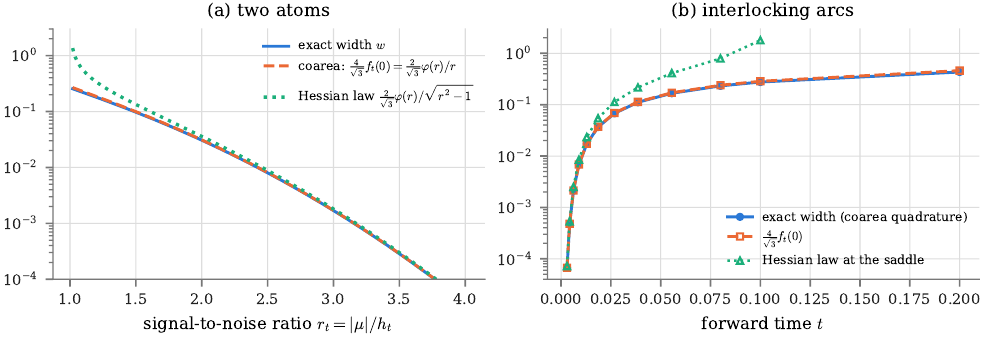}
\caption{Width of the speciation layer. (a) Two equal atoms: exact width \eqref{eq:twoatomwidth}, the coarea prediction $\frac{4}{\sqrt3}f_t(0)$ (asymptotically exact), and the saddle-curvature (Hessian) heuristic, which diverges at the speciation threshold $r_t=1$. (b) Interlocking arcs: the same quantities (exact width by coarea quadrature).}\label{fig:width}
\end{figure}

For the interlocking arcs, we compute $w(t)$ via ten-point Gauss--Legendre quadrature of the coarea integral in Theorem~\ref{thm:width}(a), projecting grid contours onto each level set by three Newton steps along \eqref{eq:gradLambda}. Point reflection $x\mapsto-x$ exchanges the classes, but since \eqref{eq:sep} fails, \eqref{eq:widthformula} is tested numerically. Figure~\ref{fig:width} compares $w(t)$ with $\frac{4}{\sqrt3}f_t(0)$ and with the saddle-curvature heuristic $\frac{2}{\sqrt3}\lambda_{\max}^{-1/2}\int_{\Gamma_t}p_t\,\dd\mathcal H^{d-1}$ of Remark~\ref{rem:widthlaw} (evaluated at the point of $\Gamma_t$ of largest density) for both two equal atoms and the interlocking arcs. On the arcs, $\frac{4}{\sqrt3}f_t(0)$ has relative error $0.4\%$ at $t=0.003$ and $3.4\%$ at $t=0.1$, whereas the heuristic overestimates $w(t)$ by $8\%$ at $t=0.003$ and by factors $1.66$ and $6.5$ at $t=0.027$ and $t=0.1$, becoming undefined once the origin turns into a maximum ($t\approx0.11$). Meanwhile, the interfaces $e^t\Gamma_t$ converge to the support bisector (Figure~\ref{fig:moons}b--c), with Hausdorff distance in $W=[-2,2]\times[-1.6,1.6]$ dropping from $0.40$ at $h_t=0.47$ to $0.023$ at $h_t=0.078$.

We also test a convexly separated example that satisfies all hypotheses of Theorems~\ref{thm:interface}(d) and~\ref{thm:width}(c). Let $\nu_1$ place masses $0.2,0.5,0.3$ at $(1,-0.5),(1.15,0.25),(1.3,0.6)$, let $\nu_2$ be its image under $y\mapsto-y$, and set $w_1=w_2=\frac12$; here \eqref{eq:sep} holds with $v=e_1$ and $g=2$, and the lower mass condition holds with $m=0$. Across $151$ equally spaced points $z\in[-1.5,1.5]$, $\max_z|\gamma^\sharp_t(z)-\gamma_0(z)|$ drops from $0.183$ at $h_t=0.6$ to $1.17\times10^{-3}$ at $h_t=0.05$, remaining in $[0.46\,h_t^2,0.63\,h_t^2]$ over all ten tested noise levels (confirming the $m=0$ rate of Theorem~\ref{thm:interface}(d)). The relative error of $\frac{4}{\sqrt3}f_t(0)$ is $0.18\%$, $0.66\%$, and $2.5\%$ at $h_t=0.2,0.4,0.8$.

\subsection{ODE basins, Bayes cells, and optimal transport}\label{subsec:num-ode}

\begin{figure}[t]
\centering
\includegraphics[width=0.64\textwidth]{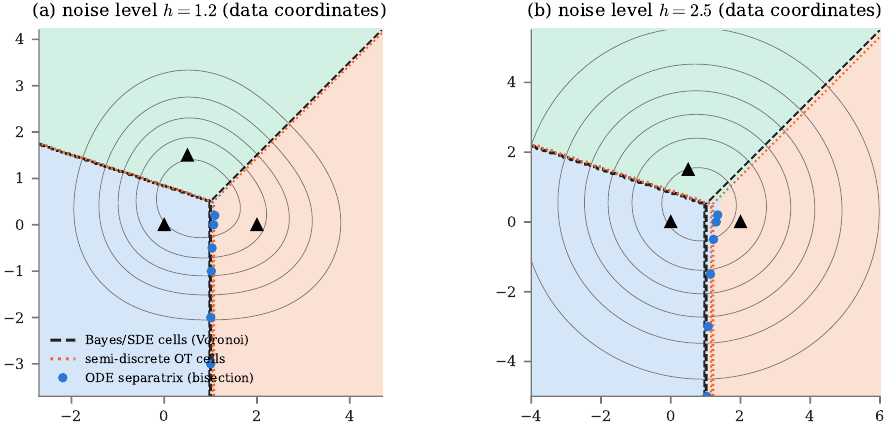}
\caption{Three equal atoms (triangles) in data coordinates. Shading: ODE basins at noise level $h$. Dashed black: SDE Bayes cells (Voronoi). Dotted orange: semi-discrete optimal-transport cells of mass $\frac13$. Dots: ODE separatrix between the two lower atoms (bisection). Thin curves: level sets of $p^h$.}\label{fig:basins}
\end{figure}

For three equal atoms at $(0,0)$, $(2,0)$, and $(\frac12,\frac32)$, we integrate the data-coordinate ODE \eqref{eq:VEflow} with RK4 in $\log h$ (step $0.002$) from $h_0\in\{1.2,2.5\}$ down to $h=0.02$ on a $301\times301$ grid, assigning each initial point $y=(y_1,y_2)$ to its nearest atom at termination. Semi-discrete optimal-transport cells of mass $\frac13$ for $p^{h_0}$ are computed via Newton's method on an entropically smoothed dual followed by hard assignment (Figure~\ref{fig:basins}). As required by Proposition~\ref{prop:odemass}, the ODE basins have masses $\frac13$ up to grid error (within $10^{-3}$), whereas the Bayes (Voronoi) cells have masses $0.326,0.343,0.331$ at $h_0=1.2$ and $0.316,0.356,0.328$ at $h_0=2.5$. The ODE partition differs from the Voronoi partition on $1.2\%$ and $2.6\%$ of the mass, and from the optimal-transport partition on $0.8\%$ and $1.3\%$ (unchanged when verified against power cells matching mass $\frac13$ to $10^{-12}$).

Theorem~\ref{thm:threeatoms} applies with $x_1=(0,0)$ opposite the longest side: the ODE partition is not a power diagram, and the basin of $x_1$ strictly contains its Voronoi cell, consistent with its Voronoi masses $0.326$ and $0.316<\frac13$. Bisection to $6\times10^{-11}$ confirms that the ODE separatrix between atoms $1$ and $2$ is curved---moving from $y_1=1.0013$ at $y_2=-3$ to $1.1086$ at $y_2=0.3$ for $h_0=1.2$, and from $1.0748$ to $1.3741$ for $h_0=2.5$---whereas the Voronoi edge is $y_1=1$ and every power-diagram edge between these atoms is vertical.

\section{Conclusion}\label{sec:conclusion}

We have treated speciation in reverse diffusion as a partition problem. For the exact score the SDE partition is explicit, because committors are posterior probabilities, while the probability-flow ODE produces a different partition with the same masses, which for three equally weighted atoms is an optimal-transport partition only in the equilateral case. The number of modes is not monotone in the noise level, but $K$ separated clusters have an explicit window with exactly $K$ modes. Under convex separation the interface is a real-analytic graph converging to the bisector of the supports, and the mass of the transition layer is given by a coarea formula. Two-sided bounds on clone agreement locate the speciation window (Proposition~\ref{prop:clonewindow}), and if an approximate score is used down to a positive noise level, with the exact posterior below it, the partition can change only inside the transition layer, up to the integrated score error (Corollary~\ref{cor:stability}). Open directions include a variational description of the ODE partition \cite{Pierret2026}, the width without symmetry, and committor interfaces for learned, non-gradient drifts, which call for the PDE methods of segregation and optimal-partition problems \cite{ContiTerraciniVerzini2005,CaffarelliLin2008,HHT2009,Boz2015Partitions,Boz2022Singular}. The Python code for the experiments of Section~\ref{sec:numerics} accompanies this preprint.

\appendix
\section{Pointwise density in the width theorem}\label{app:width}

Fix $t>0$ and assume \eqref{eq:sep}. By Theorem~\ref{thm:interface}(a), $\partial_v\Lambda_t\ge c_t>0$, so for each $z\in v^\perp$ the map $s\mapsto\Lambda_t(z+sv)$ is onto $\R$ with a smooth inverse $s(z,\lambda)$. A change of variables on each line and Fubini's theorem give the version
\[
f_t(\lambda)=\int_{v^\perp}\frac{p_t(z+s(z,\lambda)v)}{\partial_v\Lambda_t(z+s(z,\lambda)v)}\dd z
\]
of the log-odds density, and likewise for $f^{(k)}_t$. The numerator is bounded by $A\exp(-a|z|^2)$ uniformly in $s$, with $A,a>0$ depending on $t$ and the support, and the denominator is at least $c_t$; hence $f_t$ is continuous by dominated convergence, and strictly positive. Parametrizing the level surface as a graph shows that this integral is the coarea expression at every $\lambda$. The calibration and symmetry identities therefore extend by continuity to $\lambda=0$, and on $[-\log3,\log3]$ the density is bounded above and away from zero, which justifies the logarithm and the maximum in \eqref{eq:widthformula} and the pointwise substitution leading to \eqref{eq:hessbound}. Without \eqref{eq:sep} these pointwise statements can fail: in $\R^2$, for equal weights and $\nu_1,\nu_2$ uniform on $\{\pm e_1\}$ and $\{\pm e_2\}$, the log-odds in data coordinates are $\log\cosh(y_1/h^2)-\log\cosh(y_2/h^2)$, whose zero level $\{y_1=\pm y_2\}$ crosses itself at the only critical point, and $f_t(\lambda)\sim2h^4p^h(0)\log(1/|\lambda|)$ as $\lambda\to0$ ($h=h_t$), so that $f_t(0)$ is not defined.

\end{document}